\documentclass[11pt]{article}
\usepackage{benstyle}
\graphicspath{{Figures/}}
\usepackage{tcolorbox}
\newcommand{\newOmega}{Z}  

\theoremstyle{assumption}
\newtheorem{assumption}[theorem]{Assumption}

\title{An adaptive sampling and domain learning strategy for multivariate function approximation on unknown domains}

\author{Ben Adcock, Juan M. Cardenas and Nick Dexter \\ Department of Mathematics \\ Simon Fraser University \\ Canada}

\begin{document}

\newpage

\setcounter{page}{1}

\maketitle

\begin{abstract} 
Many problems arising in computational science and engineering can be described in terms of approximating a smooth function of $d$ variables, defined over an unknown \textit{domain of interest} $\Omega\subset \mathbb{R}^d$, from sample data. 
Here both the underlying dimensionality of the problem (in the case $d\gg 1$) as well as the lack of domain knowledge---with $\Omega$ potentially irregular and/or disconnected---are confounding factors for sampling-based methods. 
Na\"{i}ve approaches for such problems often lead to wasted samples and inefficient approximation schemes. For example, uniform sampling can result in upwards of 20\% wasted samples in some problems considered herein. In applications such as surrogate model construction in computational {\em uncertainty quantification} (UQ), the high cost of computing samples necessitates a more efficient sampling procedure. 
Over the last several years methods for computing such approximations from sample data have been studied in the case of irregular domains, and the advantages of computing sampling measures depending on an approximation space $P$ of $\dim(P)=N$ have been shown. 
More specifically, such approaches confer advantages such as stability and well-conditioning, with an asymptotically optimal sample complexity scaling $\mathcal{O}(N\log(N))$.
The recently-proposed {\em adaptive sampling for general domains} (ASGD) strategy is one such technique to construct these sampling measures.
The main contribution of this paper is a procedure to improve upon the ASGD approach by adaptively updating the sampling measure in the case of unknown domains. 
We achieve this by first introducing a {\em general domain adaptivity strategy} (GDAS), which computes an approximation of the function and domain of interest from the sample points. 
Second, we propose an adaptive sampling strategy, termed {\em adaptive sampling for unknown domains} (ASUD), which generates sampling measures over a domain that may not be known in advance, based on the ideas introduced in the ASGD approach. 
We then derive (weighted) {\em least squares} and {\em augmented least squares} techniques for polynomial approximation on unknown domains. 
We present numerical experiments demonstrating the efficacy of the adaptive sampling techniques with {\em least squares}-based polynomial approximation schemes. Our results show that the ASUD approach consistently achieves the same or smaller errors as uniform sampling, but using fewer, and often significantly fewer function evaluations.
\end{abstract}

%% REQUIRED
%\begin{keywords}
%high-dimensional approximation, sampling strategy, irregular domains, domain learning, surrogate model construction
%\end{keywords}

%%%%%%%%%%%%%%%%%%%%%%%%%%%%%%%%%%%%%%%%%%%%%%%%%%%%%%%%%%%%%%%%%%%%%%%
%%%----------------             Introduction        ----------------%%%
%%%%%%%%%%%%%%%%%%%%%%%%%%%%%%%%%%%%%%%%%%%%%%%%%%%%%%%%%%%%%%%%%%%%%%%

\section{Introduction}
\label{sec:intro}
In this work, we consider the problem of approximating a smooth, multivariate function $f : \Omega \rightarrow \bbR$ over a \textit{domain of interest} $\Omega \subset \bbR^{d}$ in $d \geq 1$ dimensions that is unknown \textit{a priori} and may be irregular and/or disconnected. Throughout, we assume that $\Omega$ is a subset of a known domain $D \subseteq \bbR^d$, and that for any point in $\bm{y} \in D$ we have access to a black-box that evaluates the function $f$. Given input $\bm{y} \in D$, this black box either returns the value $f(\bm{y}) \in \bbR$, or, if $f$ is undefined at $\bm{y}$, returns some form of exit flag, in which case we formally write $f(\bm{y}) = + \infty$.  Further, we assume that there is a known function  $\cQ : \bbR \cup \{ \infty\} \rightarrow \{ 0 ,1\}$ that characterizes the domain of interest $\Omega$ through evaluations of $f$. In other words, we may write
\be{
\label{Omega_Q}
\Omega = \{ \bm{y} \in D :  \cQ(f(\bm{y})) = 1 \},
}
where, in particular, $\cQ(+\infty) = 0$. Typical examples we consider include the case
\bes{
\Omega = \{ \bm{y} \in D : f(\bm{y}) < \infty \},
}
i.e.\ where $\Omega$ is precisely the domain within $D$ in which $f$ is defined, or
\be{
\label{Omega_examples}
\Omega = \{ \bm{y} \in D : 0 \leq f(\bm{y}) < \infty \},
}
in which case $\Omega$ is the domain of interest in which $f$ is nonnegative. Note that both of these are examples of the general form
\be{
\label{Omega_examples_2}
 \Omega = \{ \bm{y} \in D : a \leq f(\bm{y}) < b\},
}
for constants $- \infty < a < b \leq \infty$. Observe that in this example, the function $\cQ$ can be taken as the indicator function $\bbI_{[a,b)}$ of the interval $[a,b)$. In general, $\cQ$ can be chosen as the indicator function $\bbI_{f(\Omega)}$ of the set $f(\Omega)$. However, it is convenient for the purposes of this paper to define the set $\Omega$ in terms of a function $\cQ$, as in \R{Omega_Q}, rather than the other way around.

Our work is motivated by applications in \textit{Uncertainty Quantification (UQ)}. For this reason, we focus on three main challenges: first, the \textit{curse of dimensionality}, since $d$ is often much larger than one in practice; second, computing pointwise evaluations of $f$ via the black box is expensive, and; third, the domain of interest $\Omega$ is unknown in advance, only accessible through evaluating $f$ (see \eqref{Omega_Q}) and is generally irregular.

\subsection{Motivations} 

Our primary motivation is surrogate model construction in UQ, in which case $f$ is some unknown \textit{quantity of interest}, $\bm{y}$ represents the parameters in the model and the black box is typically some large computer code that evaluates $f$, e.g., a numerical PDE solver. Such problems are often high dimensional, since models usually involve many parameters, and samples are expensive to compute, since they typically involve expensive computer simulations, (see, e.g., \cite{Smith2013,Sullivan2015} for introductions to UQ, as well as \cite{adcock2021sparse}). 

It is standard in surrogate model problems to assume that $D$ is a tensor-product domain, e.g.\ the unit hypercube $[-1,1]^d$, and then to perform \textit{Monte Carlo sampling} over $D$. That is, the sample points $\bm{y}_1,\bm{y}_2,\ldots$ are drawn randomly and independently from an underlying probability distribution on $D$, such as the uniform or Chebyshev distribution. 
However, in practice the situation often arises where $f$ represents some physical quantity (mass, pressure, and so forth) which may be known to be nonnegative (e.g.\ $\Omega$ as in \R{Omega_examples}) or bounded between finite minimum and maximum values (e.g.\ $\Omega$ as in \R{Omega_examples_2}), in which case, the model has no physical interpretation outside of $\Omega$.  In other cases, the model may simply not be well defined over the whole of $D$, e.g.\ in the context of a PDE model, if certain parameter values lead to an ill-posed PDE system. In general, though, there may be no straightforward way to determine $\Omega$ without first evaluating $f$. This makes Monte Carlo sampling wasteful, since any sample that returns a value $\cQ(f(\bm{y}_i)) = 0$ (in the general setting \eqref{Omega_Q}) is simply rejected. Indeed, one expects roughly $\mathrm{meas}(D\backslash \Omega)/\mathrm{meas}(D)$ of the evaluations $f(\bm{y}_1),f(\bm{y}_2),\ldots$ to return such a value (here $\mathrm{meas}$ is the measure of a set with respect to the probability measure) thus wasting a constant proportion of the samples. In practical situations, it is not uncommon for 15\% to 20\% of the samples to be wasted in this way \cite{MullerDay2018,SargsyanEtAl2014}. 

Note that this problem is related to the problem of \textit{dependent random variables} in surrogate model construction \cite{ernst_convergence_2012,jakeman2019polynomial, le2010spectral,soize_physical_2004}. It is also related to the problem of \textit{hidden constraints} in surrogate optimization \cite{AudetEtAl2020,KyEtAl2016,LeeEtAl2011,MullerDay2018}.

In this paper we introduce a new procedure for approximating a function $f$ efficiently over an unknown domain of interest. Our approach is based on an adaptive sampling procedure, termed \textit{Adaptive Sampling for Unknown Domains (ASUD)}, that improves on function approximation strategies based on Monte Carlo sampling. It does so by iteratively learning an approximation to the domain of interest $\Omega$, and then using this information to adaptively define new sampling distributions from which subsequent sample points are drawn.

\subsection{Adaptive Sampling for General Domains (ASGD)} 

Our approach is based on recent work on adaptive sampling for irregular domains. In \cite{AdcockCardenas2020}, the first two authors introduced a method for function approximation -- that we henceforth refer to as \textit{Adaptive Sampling for General Domains (ASGD)} -- over arbitrary, but known domains $\Omega$ (similar approaches have also been developed in \cite{dolbeault2020optimal,MiglioratiIrregular}). The essence of this method is a (weighted) least-squares procedure in an arbitrary finite-dimensional space $P$, typically a polynomial space, of dimension $\dim(P)=N$. It first replaces $\Omega$ by a fine grid of $K \gg N$ points. Next, it defines certain discrete sampling measures over the grid, related to the \textit{Christoffel function} of $P$, from which the samples are drawn randomly and independently. As was shown in \cite{AdcockCardenas2020}, this method has provably near-optimal sample complexity, with the number of samples $M$ required for a quasi-best approximation $\tilde{f} \in P$ to $f$ scaling like $N \log(N)$. To make this procedure adaptive, one first considers a sequence of nested subspaces $P_1\subset P_2\subset\ldots$ of dimensions $\dim(P_k) = N_k$, $k = 1,2,\ldots$. Then, following an approach introduced in \cite{migliorati2019adaptive}, one generates a sequence of sampling measures so that at the $(k+1)$th step a total of $M_{k+1} - M_{k}$ new samples are drawn and combined with the existing $M_{k}$ samples, giving a total of $M_{k+1}$ samples. As shown in \cite{AdcockCardenas2020}, if the subspaces are nested and the measures defined in a suitable way, a sequence of quasi-best approximations $\tilde{f}_k \in P_k$ is obtained from a near-optimal sample complexity, i.e.\ $M_k \asymp N_k \log(N_k)$ for each $k$. Note that alternative adaptive schemes are also possible; see, e.g., \cite{arras2019sequential}.

\subsection{Contributions}
\label{sec:con}
 
In this work, we extend the ASGD method to tackle the significantly more challenging scenario where the domain of interest $\Omega$ is unknown in advance. Our main contribution is the introduction of the aforementioned ASUD  procedure. This method is in turn based on a general strategy for function approximation and domain learning, which we term \textit{General Domain Adaptivity Strategy (GDAS)}. Similar to previous works (see, e.g.,\cite{dolbeault2020optimal,AdcockCardenas2020,MiglioratiIrregular}) the approximation is computed via a (weighted) least-squares procedure. However, we consider two different formulations of this procedure, termed ASUD-LS (ASUD-\textit{least squares}) and ASUD-ALS (ASUD-\textit{augmented least squares}). We discuss the relative merits of each procedure. Specifically, the latter can offer better domain learning in practice, at the price of worse function approximation when the function is badly behaved (e.g.\ singular) outside of the domain of interest.
 
We present a series of numerical experiments to demonstrate the benefits of ASUD on different problems. In these experiments, we observe significant benefits of ASUD over \textit{Monte Carlo (MC)} sampling, both in terms of approximating the function and learning the domain. Even in the best cases, MC sampling requires at least 50\% more samples to achieve a similar error to ASUD. In other cases, it may also fail to achieve the same accuracy. On the other hand, while MC sampling wastes a significant proportion of the samples, the ASUD procedure is asymptotically optimal in terms of its samples. That is to say, the proportion of rejected samples tends to zero as the number of iterations increases. In fact, the performance of ASUD is very similar to ASGD, which requires \textit{a priori} domain knowledge. The GitHub repository can be found in \url{https://www.github.com/JMcardenas/Adaptive-sampling-and-domain-learning}.

\subsection{Outline} 

The outline of the remainder of this paper is as follows. We commence in \S \ref{sec:GDAS} by introducing GDAS. Next, in \S \ref{sec:ASUDmain} we first review ASGD and then introduce ASUD. In \S \ref{sec:other_decoders} we then introduce the two approximation methods studied: (weighted) least squares and augmented (weighted) least squares. In \S \ref{sec:num_exp} we present a series of numerical experiments comparing the various methods. Then in \S \ref{sec:theoretical} we discuss a number of theoretical considerations. Finally, we conclude in \S \ref{sec:conclusions} with some open problems and topics for future work.
 
\section{General Domain Adaptivity Strategy (GDAS)}\label{sec:GDAS}
In this section, we describe the GDAS procedure on which the methods developed later are based.

\subsection{Setup}
As noted above, we consider a domain $D \subseteq \bbR^d$, an unknown black box function $f : D \rightarrow \bbR \cup \{ + \infty\}$  and a known function $\cQ : \bbR \cup \{ + \infty\} \rightarrow \{ 0 ,1\}$ with $\cQ(+\infty) = 0$ that characterizes an unknown domain of interest $\Omega \subseteq D$ as
\be{\label{def:Omega}
\Omega = \{ \bm{y} \in D :  \cQ(f(\bm{y})) = 1 \}.
}
Following the approach of \cite{AdcockCardenas2020,MiglioratiIrregular,dolbeault2020optimal}, our first step is to discretize the domain $D$. As shown in these works, discretizing the domain allows for the construction of optimal sampling measures in the setting of possibly irregular domains. To this end, we assume that there is a finite grid of points $Z \subset D$ and a discrete probability measure $\tau$ supported on $Z$. Normally, $\tau$ is taken to be the discrete uniform measure over $Z$. Throughout the paper, we assume samples $\bm{y}_1,\bm{y}_2,\ldots$ are drawn randomly according to certain distributions that are supported on $Z$. We write $\tilde{\mu}_i$ for the measure supported on $Z$ from which the $i$th sample $\bm{y}_i$ is drawn.

The grid $Z$ and measure $\tau$ serve two purposes. First, $Z$ is used to perform the key computations: namely, adaptively generating the discrete sampling measures in the algorithm. 
Second, $\tau$ -- or more precisely its normalized restriction to $\Omega$, which we denote by $\tilde{\tau}$ -- is the measure with respect to which we evaluate the error of the approximation. To be precise, if $\tilde{f}$ is an approximation to $f$, we compute the error
\bes{
E(f)=\frac{\nmu{f-\tilde{f}}_{L^2(\Omega,\tilde{\tau})}}{\nm{f}_{L^2(\Omega,\tilde{\tau})}},\quad
\mathrm{d}\tilde{\tau}(\bm{y}):=\frac{\bbI_{\Omega}(\bm{y})}{\int_\Omega d\tau(\bm{y})}\mathrm{d}\tau(\bm{y}).
}
 
Our objective is to approximate $f$ over the unknown domain $\Omega$ as accurately as possible (with respect to this error measure) from as few samples as possible. Hence, we assume throughout that $Z_{\Omega} : = Z \cap \Omega$ is sufficiently fine so as to represent $\Omega$ well. Loosely speaking, by this we mean that if an approximation $\tilde{f}$ approximates $f$ well on $Z_{\Omega}$ then it also approximates $f$ well on $\Omega$ itself. This raises the question of how fine the grid $Z$ should be in practice. We will not discuss this issue, although we note that theoretical estimates are available in many cases -- in particular, when the grid is generated via Monte Carlo sampling. See 
\cite{adcock2020approximating,dolbeault2020optimal} for further discussion on this topic.

\subsection{GDAS}

Given $Z$ and $\tau$, we first define \textit{sampling numbers} $0 = M_0 < M_1 <M_2<\ldots$, where $M_l$ is the number of number of samples of $f$ from $\Omega$ used in the $l$th step. We also set $Z_0=Z$. The $l = 1$ step of the GDAS strategy now proceeds as follows: 

\begin{itemize}
\item[(a)] We first construct sampling measures $\mu_{1},\ldots,\mu_{M_1}$ supported on $Z_0$.
\item[(b)] Next, for each $i = 1,\ldots,M_1$, we draw a point $\bm{y}_i$ randomly according to $\mu_{i}$ and evaluate $f(\bm{y})$. If $\cQ(f(\bm{y})) = 1$ then we accept $\bm{y}$ and write $\bm{y}_i = \bm{y}$. Otherwise, we reject it, draw a new point $\bm{y}$ and repeat until we obtain $\cQ(f(\bm{y})) = 1$. Having done this for each $i = 1,\ldots,M_1$ and obtained the sample points $\cS_1 = \{\bm{y}_i \}^{M_1}_{i=1}$, we write $F_1 \geq M_1$ for the total number of evaluations of $f$ used and let $\cR_1 = \{ \bm{u}_i \}^{T_1}_{i=1}$ be the rejected points, where $T_1 = F_1 -M_1$. 
 
\item[(c)] Next, we compute an approximation $\tilde{f}_1$ to $f$ from the values $\{ f(\bm{y}_i) \}^{M_1}_{i=1}$, and potentially also using the rejected values $\{ f(\bm{u}_i) \}^{T_1}_{i=1}$ (if finite).
\item[(d)] We then compute the domain approximation
\bes{
Z_1 = \{ \bm{z} \in Z : \cQ(\tilde{f}_1(\bm{z})) = 1 \} \cup \cS_1 \backslash \cR_1 
}
using the function $\cQ$ applied to the approximation $\tilde{f}_1$. We explicitly include the sample points $\cS_1 = \{ \bm{y}_i \}^{M_1}_{i=1}$ and exclude the rejected points $\cR_1 = \{ \bm{u}_i \}^{T_1}_{i=1}$.
\end{itemize}
Note that one may choose the initial measures simply as the discrete measure over $Z$, i.e. $\mu_1 = \ldots = \mu_{M_1} = \tau$. This is common in surrogate model construction. However, as we explain in Section \ref{sec:ASUD} one may also consider other approaches.

Having completed the $l = 1$ step, the $l$th step of GDAS, $l \geq 2$, proceeds as follows:
\begin{itemize}
\item[(a)] We use the domain approximation $Z_{l-1}$ to construct sampling measures $\mu_{M_{l-1}+1},\ldots,\mu_{M_{l}}$ supported on $Z_{l-1}$.
\item[(b)] For each $i = M_{l-1}+1,\ldots,M_l$, we proceed as in (a) above, drawing $\bm{y}_i$ using a combination of $\mu_i$ and rejection sampling according to $\cQ(f(\bm{y}))$. We write $G_{l} \geq M_l - M_{l-1}$ for the total number of function evaluations in doing this, so that $F_l = F_{l-1} + G_l$ is the total number of function evaluations used up to and including the $l$th step. We also let $\cS_l = \{ \bm{y}_i \}^{M_l}_{i=1}$ be the set of all sample points and $\cR_l = \{ \bm{u}_i \}^{T_l}_{i=1}$ be the set of all rejected points, where $T_l = F_l -M_l$.
\item[(c)] We compute an approximation $\tilde{f}_l$ to $f$ from the values $\{ f(\bm{y}_i) \}^{M_l}_{i=1}$, and potentially also using the rejected values $\{ f(\bm{u}_i) \}^{T_l}_{i=1}$ (if these values are finite).
\item[(d)] We compute the domain approximation
\bes{
Z_l = \{ \bm{z} \in Z : \cQ(\tilde{f}_l(\bm{z})) = 1 \} \cup \cS_l \setminus \cR_l. 
}
\end{itemize}

\begin{figure}
\begin{center}
\includegraphics[scale=0.5]{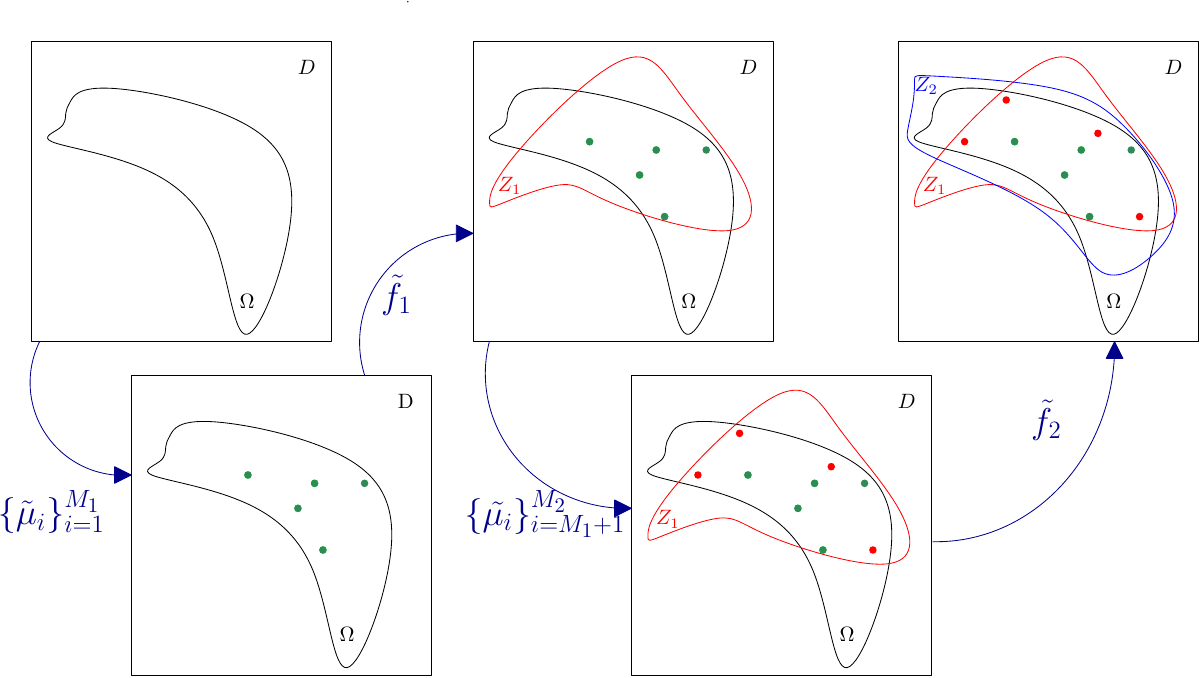}
\caption{A visualization of two steps of GDAS.}
\label{fig:GDAS_iter1}
\end{center}
\end{figure}

This procedure is summarized as Method 1 and visualized in Fig.\ \ref{fig:GDAS_iter1}.

\begin{tcolorbox}[floatplacement=t!,float,title=Method 1. General Domain Adaptivity Strategy (GDAS)]
\noindent \textbf{Inputs:} Finite grid $Z$, probability measure $\tau$ over $Z$, sampling numbers $0 = M_0 < M_1 < M_2 < \ldots $, function $\cQ$ as \ref{def:Omega}.
\\
\noindent \textbf{Initialize:} Set $\newOmega_0 = Z$, $F_0 = 0$, $\cS_1 = \cR_1 = \emptyset$.\\ 
\\
\noindent \textbf{for} $l=1,2,\ldots$ \textbf{do}\\
\textbf{Stage (a)} Construct measures $\{ \mu_i \}^{M_{l}}_{i= M_{l-1} + 1}$ over $\newOmega_{l-1}$.
\\
\textbf{Stage (b)} \\
\noindent \textbf{for} $i=M_{l-1}+1,\ldots,M_l$ \textbf{do}\\ 
Draw $\bm{y}$ randomly and independently from $\mu_i$. If $\cQ(f(\bm{y})) = 0$ then add $\bm{y}$ to $\cR_l$. Repeat until $\cQ(f(\bm{y}))=1$ and then add $\bm{y}$ to $\cS_l$.
\\
\noindent \textbf{end for}
\\
\noindent \textbf{return} Sample points $\cS_l = \{ \bm{y}_i \}^{M_l}_{i=1}$ and rejected points $\cR_l = \{ \bm{u}_i \}^{T_l}_{i=1}$, total number of function evaluations $F_l = F_{l-1} + G_l$, where $G_l$ is the number of evaluations of $f$ used in this stage, and $T_l = F_l - M_l$.
\\
\textbf{Stage (c)} Compute an approximation $\tilde{f}_{l}$ to $f$ from the data $\{f(\bm{y}_i)\}^{M_l}_{i=1}$, and potentially also using the rejected values $\{ f(\bm{u}_i) \}^{T_l}_{i=1}$.
\\
\textbf{Stage (d)}  Compute $\newOmega_l = \{ \bm{z} \in Z : \cQ(\tilde{f}_l(\bm{z})) = 1 \} \cup \cS_l \setminus \cR_l $. 
\\
Update $\cS_{l+1} = \cS_l$ and $\cR_{l+1} = \cR_l$
\\
\noindent \textbf{end for}
\\
\textbf{Outputs:} The sequences of function approximations $\{\tilde{f}_l\}_{l\geq 1}$ and domain approximations $\{\newOmega_l\}_{l\geq 1}$.
\end{tcolorbox}  

\subsection{Discussion}
\label{sec:setup_discussion}
Some remarks are in order. First, notice that this procedure produces sequences of approximations $\{ \tilde{f}_l \}$ and $\{Z_l\}$ to both the function $f$ \textit{and} the domain $\Omega$. We refer to the latter as \textit{domain learning}. It is also \textit{adaptive}, in that the domain estimate $Z_{l-1}$ is used in stage (a) to construct the new sampling measures. 

Second, observe that the $l = 1$ step of GDAS is, in the case that the initial measures are taken as $\mu_1 = \ldots = \mu_{M_1} = \tau$, equivalent to standard Monte Carlo sampling for surrogate model construction (albeit with a discrete measure). The samples $\bm{y}_1,\ldots,\bm{y}_{M_1}$ are drawn according to the initial measure $\tau$, combined with rejection sampling based on the value of $f(\bm{y})$. A substantial fraction of samples may therefore be rejected in this step. Or, in other words, the \textit{rejection rate} $(F_1 - M_1)/F_1$ may be large. However, by adaptively learning $Z_l$, we hope that the overall rejection rate
\bes{
\frac{F_l - M_l}{F_l},
}
decreases with increasing $l$. Ideally, it should decrease to zero as $l \rightarrow \infty$, i.e. asympototically, a vanishing proportion of samples are wasted. In this case, we say that the method is \textit{asymptotically optimal} in terms of its samples. Later, in our numerical examples, we show that this situation does indeed occur in practice.  

Third, notice that the point $\bm{y}_i$ is computed by repeatedly drawing samples independently from $\mu_i$ and rejecting them if they fall outside $\Omega$. This means that $\bm{y}_i \sim \tilde{\mu}_i$, where $\tilde{\mu}_i$ is the (normalized) restriction of $\mu_i$ to $\Omega$, i.e.
\be{
\label{tilde_mu_i_def}
\D \tilde{\mu}_i(\bm{y}) = \frac{\bbI_{\Omega}(\bm{y}) \D \mu_i(\bm{y})}{\int_{\Omega} \D \mu_i(\bm{y}) }.
}
Since, by construction, the measure $\mu_i$ is supported on $Z_l$ for $i = M_{l-1}+1,\ldots,M_l$, the measure $\tilde{\mu}_i$ is supported on $Z_l \cap \Omega$.

Fourth, notice that the grid $Z_l$ computed in (d) is constructed as the union of the points from $Z$ where $\tilde{f}_l$ returns the right value via $\cQ$ and the sample points $\cS_l = \{ \bm{y}_i \}^{M_l}_{i=1}$. The reason for defining $Z_l$ in this way is that the approximation $\tilde{f}_l$ will generally not be an interpolant of $f$ at the sample points. Hence one may have $\cQ(\tilde{f}_l(\bm{y}_i)) = 0$ for some $i$, even though $\cQ(f(\bm{y}_i)) = 1$ by construction. Therefore, to give the best estimate of $\Omega$, we include both sets of points in the definition of $Z_l$. In order to enhance the estimate of $\Omega$ further, the rejected points $\cR_l = \{ \bm{u}_i \}^{T_l}_{i=1}$ are also removed from $Z_l$, since we know they must lie outside of $\Omega$.

Finally, we note that the GDAS procedure does not specify how the approximations $\tilde{f}_l$ are computed, nor how the sampling measures $\mu_l$ are defined. This is the concern of the next two sections.

\section{Adaptive Sampling for Unknown Domains (ASUD)}\label{sec:ASUDmain}
In this section, we describe how to generate the sampling measures, leading to the ASUD method. We divide this section into two parts. First, we recap how the sampling measures can be constructed when the domain is known. This is based on the ASGD method introduced in \cite{AdcockCardenas2020,MiglioratiIrregular}. With that in hand, we then introduce ASUD.

\subsection{ASGD}\label{alg:asgd}

Suppose that $\Omega$ is known, $Z = \{ \bm{z}_i \}^{K}_{i=1} \subset \Omega$
is a fine grid that represents $\Omega$ and $\tau$ is a probability measure that is supported on $Z$. Notice that since $\Omega$ is known in this case, we take $Z$ as a grid over $\Omega$ as opposed to $D$.

The basic idea behind ASGD is the following. Let $P \subset L^2(\Omega,\tau)$ be a finite-dimensional subspace of dimension $\dim(P) = N$ in which we expect that $f$ is well approximated. For example, if $f$ is smooth we may take $P$ to be a subspace of multivariate polynomials of a given maximum degree. We also assume that $P$ satisfies the following assumption:
\be{
\label{P_assume}
\text{For any $\bm{y}\in \mathrm{supp}(\tau)$ there exists a $p\in P$ with $p(\bm{y})\neq 0$.}
}
Note that this assumption trivially holds whenever $P$ contains the constant function. This will always be the case in our numerical examples.

The goal of ASGD is to use the information about $\Omega$ and $P$ to devise sampling measures $\mu_1,\ldots,\mu_M$ -- from which the $M$ sample points $\bm{y}_1,\ldots,\bm{y}_M$ are to be drawn randomly  -- for which a log-linear sample complexity $M \asymp N \log(N)$ is sufficient to stably represent elements of $P$ through their sample values. To be precise, we mean that for any $0 < \delta < 1$ there is a constant $c_{\delta} >0$ such that if
\bes{
M \geq c_{\delta} N \log(N),
}
 then, with high probability,
when $M$ samples are drawn independently with $\bm{y}_i \sim \mu_i$ for $i = 1,\ldots,M$ one has the equivalence
\be{
\label{L2_equiv}
(1-\delta) \nm{p}^2_{L^2(\Omega,\tau)} \leq \frac1M \sum^{M}_{i=1} w(\bm{y}_i) |p(\bm{y}_i) |^2 \leq (1+\delta) \nm{p}^2_{L^2(\Omega,\tau)} ,\quad \forall p \in P,
}
over the subspace $P$ between the $L^2(\Omega,\tau)$-norm and a certain weighted discrete norm defined in terms of a positive weight function $w$.

\rem{
\label{r:samplecomplexitystable_LS}
A particular consequence of this equivalence is that, if an approximation $\tilde{f}$ is computed as the (weighted) least-squares fit from the subspace $P$ based on the sample values $\{ f(\bm{y}_i) \}^{M}_{i=1}$, then $\tilde{f}$ is a quasi-best approximation to $f$ from $P$.  Specifically, if \eqref{L2_equiv} holds then the weighted least-squares approximation 
\be{
\tilde{f} = \argmin{p\in P}\left\lbrace \frac{1}{M}\sum_{i=1}^M w(\bm{y}_i) |f(\bm{y}_i)-p(\bm{y}_i))|^2\right\rbrace,
}
(which is unique)
of a function $f\in L^2(\Omega,\tau)$ satisfies
\bes{
\nmu{f-\tilde{f}}_{L^2(\Omega,\tau)}\leq  \inf_{p \in P} \left \{ \nm{f-p}_{L^2(\Omega,\tau)} + c'_{\delta} \nm{f-p}_{\mathrm{disc}} \right \},
} 
for some $c'_{\delta} > 0$ depending on $\delta$, where $\nm{g}_{\mathrm{disc}} = \sqrt{M^{-1} \sum^{M}_{i=1} w(\bm{y}_i) | g(\bm{y}_i) |^2 }$. 
}

The construction of such sampling measures is achieved via the \textit{Christoffel function} of $P$. If $\{ \phi_i \}^{N}_{i=1} \subset P$ is any orthonormal basis of $P$ in $H := L^2(\Omega,\tau)$, we define the normalized, reciprocal Christoffel function {$\mathcal{K}(P,H)$} as
\be{
\label{Christoffel_fun}
{\cK(P,H)(\bm{y})} = \frac1N \sum^{N}_{i=1} |\phi_i(\bm{y}) |^2,
}
(the function $1/\sum_{i=1}^N|\phi_i(\bm{y})|^2$ is the Christoffel function of $P$ in $H$; see \cite{NevaiFreud1986}).
Note that {$\cK(P,H)$} is strictly positive for each $\bm{y} \in \mathrm{supp}(\tau)$ if and only if \eqref{P_assume} holds.

The idea of using the Christoffel function to construct the sampling measures was considered in \cite{hampton2015coherence} for total degree polynomial spaces, then later in \cite{MiglioratiCohenOptimal} for arbitrary spaces. See \cite[Sec.\ 5.5]{adcock2021sparse} for an overview. As was shown in \cite{AdcockCardenas2020,MiglioratiCohenOptimal}, a suitable choice of sampling measures is any collection $\{ \mu_i \}^{M}_{i=1}$ that satisfies
\be{
\frac{1}{M}\sum_{i=1}^{M}\D\mu_i(\bm{y}) = {\cK(P,H)(\bm{y})} \D \tau(\bm{y}) = 
\frac{1}{N}\sum_{i=1}^{N}\left|\phi_i(\bm{y})\right|^2\D\tau(\bm{y}),\quad \bm{y} \in \Omega.\label{eq:Christoffel_cond}
}
When chosen in this way, \eqref{L2_equiv} holds with the weight function given by $w(\bm{y}) = \left ( {\cK(P,H)(\bm{y})} \right )^{-1}$.
In practice, there are several ways to choose measures $\mu_i$ that satisfy \eqref{eq:Christoffel_cond}. One option is simply to set
\bes{
\D \mu_i(\bm{y}) = {\cK(P,H)(\bm{y})} \D \tau(\bm{y}),\quad i = 1,\ldots,M.
}
However, this approach is not well suited to hierarchical schemes, in which, rather than a fixed subspace $P$ one seeks a sequence of approximations in a nested family of subspaces $P_1 \subseteq P_2 \subseteq \cdots$. See \cite{arras2019sequential,migliorati2019adaptive} and \cite[Sec.\ 5.5]{adcock2021sparse} for further discussion.

To develop a hierarchical procedure, we proceed as follows. First, we let $M = k N$ for some sampling ratio $k \in \bbN$. Then we choose the sampling measures $\mu_i$ as
\be{\label{measure:mu_i}
\D\mu_i(\bm{y}) = \left| \phi_j(\bm{y})\right|^2\D\tau(\bm{y}),\quad\bm{y}\in\Omega,\ (j-1)k<i\leq jk,\ j=1,\ldots,N.
}
In other words, the first $k$ points are drawn from a measure weighted by $|\phi_1(\bm{y})|^2$, the next $k$ points from a measure weighted by $|\phi_2(\bm{y})|^2$, and so forth.
The choice \eqref{measure:mu_i} readily leads to an hierarchical scheme. Indeed, suppose the space $P$ is enriched to a new space $\tilde{P}$ of dimension $\dim(\tilde{P}) : = \tilde{N} \geq N$ with $P \subseteq \tilde{P}$. Let
\bes{
\{ \phi_1,\ldots,\phi_N , \phi_{N+1} , \ldots, \phi_{\tilde{N} } \} \subset L^2(\Omega,\tau),
}
be an orthonormal of $\tilde{P}$, where $\{\phi_i\}_{i=1}^N$ is the original orthonormal basis of $P$. Consider a new sampling ratio $\tilde{k}\in\mathbb{N}$, $k\leq \tilde{k}$ and define $\tilde{M}=\tilde{k}\tilde{N}$. Then we can retain the original $M$ measures $\mu_1,\ldots,\mu_M$ and define new measures $\mu_{M+1},\ldots,\mu_{\tilde{M}}$ by
\eas{
\D\mu_i(\bm{y}) &= |\phi_j|^2\D\tau(\bm{y}), \quad
\bm{y}\in\Omega,\ (j-1)(\tilde{k}-k)< i - M \leq j(\tilde{k}-k),\ 
j=1,\ldots,N,\\
\D\mu_i(\bm{y}) &= |\phi_j|^2\D\tau(\bm{y}), \quad
\bm{y}\in\Omega,\ (j-1)\tilde{k}<i\leq j\tilde{k},\ j=N+1,\ldots,\tilde{N}.
} 
It is readily seen that the augmented set of measures $\mu_1,\ldots,\mu_{\tilde{M}}$ satisfies the condition \eqref{eq:Christoffel_cond} for the augmented space $\tilde{P}$. However, since the first $M$ measures remain unchanged, this means we can recycle the previously-drawn samples $\bm{y}_1,\ldots,\bm{y}_M$ when constructing the subsequent approximation over $\tilde{P}$, i.e.\ we only need to draw $\tilde{M} - M$ new samples $\bm{y}_{M + 1},\ldots,\bm{y}_{\tilde M}$ according to the new measures.

\rem{
The above discussion assumes an orthonormal basis for $P$ with respect to $\tau$. This is needed, in particular, to define the sampling measures in \eqref{measure:mu_i}. In practice, such a basis may not be given in advance, but it can be readily computed via QR decomposition. This is, in fact, the reason for defining $\tau$ as a measure supported on a finite grid as opposed to the whole of $\Omega$. We now describe the construction of such a basis and how to sample from the corresponding measures \eqref{measure:mu_i}.

Since $\tau$ is a measure defined over $Z$ we can write
\bes{
\D\tau(\bm{y})=\sum_{i=1}^K \tau_i \delta(\bm{y}-\bm{z}_i),
}
where $\tau_i = \bbP(\bm{y} = \bm{z}_i)$ for $\bm{y} \sim \tau$. Now let $\{ \psi_i \}^{N}_{i=1}$ be a nonorthogonal basis for $P$,
\bes{
\bm{B} = \left\lbrace \sqrt{\tau_i} \psi_j(\bm{z}_i)\right\rbrace_{i,j=1}^{K,N}\in \mathbb{C}^{K\times N},
}
and suppose that $\bm{B}$ has QR decomposition $\bm{B}=\bm{Q}\bm{R}$, where $\bm{Q}\in\bbC^{K\times N}$ and $\bm{R}\in\bbC^{N \times N}$. Write $\bm{Q}=\{q_{ij}\}_{i,j=1}^{K,N}$. Then it follows straightforwardly that 
\bes{
\phi_i(\bm{y}) = \sum^{i}_{j=1} (\bm{R}^{-\top})_{ij} \psi_j(\bm{y}),\quad i = 1,\ldots,N.
}  
Thus, the orthonormal basis functions $\phi_i$ are easily computed from the $\psi_j$'s. Moreover, if $M = k N$ and $(j-1) k < i \leq j k$, then the measure $\mu_i$ defined in \eqref{measure:mu_i} is given by
\bes{
\D\mu_i(\bm{y})
=|\phi_j(\bm{y})|^2\D\tau(\bm{y}) 
= \sum_{k=1}^K \tau_k |\phi_j(\bm{z}_k)|^2\delta(\bm{y}-\bm{z}_k)\D\bm{y} = \sum^{K}_{k=1} | q_{kj} |^2 \delta(\bm{y}-\bm{z}_k) \D \bm{y}.
}
Hence, $\mu_i$ is the discrete measure with $\bm{y} \sim \mu_i$ if $\bbP(\bm{y} = \bm{z}_k) = |q_{kj}|^2$ for $k = 1,\ldots,K$. Thefore, sampling from $\mu_i$ is equivalent to drawing an integer $l$ randomly from the set $\{1,\ldots,K\}$ according to the distribution $\{ |q_{kj}|^2 \}^{K}_{k=1}$ and then setting $\bm{y} = \bm{z}_l$.  
}
 
\subsection{ASUD}\label{sec:ASUD}
We now return to the main problem considered in this paper -- namely, where $\Omega$ is unknown -- and introduce the ASUD method. ASUD is an instance of GDAS in which we use the ideas from ASGD to construct the sampling measures. To this end, we consider the following setup: 

\begin{itemize}
\item[•] Let $0 = k_0 < k_1 \leq k_2 \leq \ldots $ be a sequence of sampling ratios and $0 = N_0 < N_1 \leq N_2 \leq N_3 \leq \ldots$ be a sequence of integers  such that $M_l = k_l N_l$ for $l = 0,1,2,\ldots$.
\item[•] Let $P_1 \subseteq P_2 \subseteq P_3 \subseteq \cdots \subset L^2(D,\tau)$ be a nested sequence of finite-dimensional subspaces of dimensions $\dim(P_l) = N_l$ and $\psi_1,\psi_2,\ldots$ be functions such that $P_l = \spn \{ \psi_{i} \}^{N_l}_{i=1}$ for $l = 1,2,\ldots$.  
\end{itemize} 
Note that the subspaces $P_l$ may either be defined beforehand or generated adaptively during GDAS. The idea behind the second part of this setup is that the $P_l$ are subspaces in which we expect $f$ to be well approximated, i.e.\ the best approximation error decreases rapidly in $l$. In this paper, we consider the case where the subspaces are defined \textit{a priori} beforehand. However, it is also possible to generate them adaptively. For example, in the case of polynomial approximation, one often does this in a greedy manner. Here, $P_l$ is a polynomial subspace spanned by polynomials with multi-indices in a given set $S_l \subset \bbN^d_0$. At step $l$, the computed approximation $\tilde{f}_l \in P_l$ is used to estimate the polynomial coefficients for indices lying on the so-called \textit{reduced margin} of $S_l$. Then, the indices corresponding to the largest coefficients in magnitude are used to construct the next index set $S_{l+1}$ and corresponding polynomial subspace $P_{l+1}$. See, e.g., \cite{cohen2018multivariate,migliorati2015adaptive,migliorati2019adaptive} for further information.
%typically in a greedy manner (see, e.g., \cite{cohen2018multivariate,migliorati2015adaptive,migliorati2019adaptive}).

With this setup in hand, we now describe how ASUD constructs the sampling measures $\mu_1,\mu_2,\ldots$. Consider the first step $l = 1$ of GDAS, in which we first need to define the measures $\mu_1,\ldots,\mu_{M_1}$. We do this via the same approach as ASGD. For convenience, we now also write $\Omega_0 = D$ for the first estimate of the domain $\Omega$, $Z_0 = Z$ for the first estimate of the discrete domain $\Omega \cap Z$ and $\tau_0 = \tau$. Let $P_1$ be the first subspace. We construct an orthonormal basis $\{ \phi^{(1)}_{j} \}^{N_1}_{j=1}$ for $P_1$ over $Z_0$, so that
\bes{
P_1 = \spn \{ \phi^{(1)}_{1},\ldots,\phi^{(1)}_{N_1} \} \subset L^2(Z_0,\tau_0) \equiv L^2(D,\tau).
} 
Then we define the measures exactly as in ASGD, i.e.
\bes{
\D\mu_{i}(\bm{y}) = | \phi^{(1)}_{j}(\bm{y}) |^2 \D \tau_0(\bm{y}),\quad\bm{y}\in D,\ (j-1)k_1<i\leq jk_1,\ j=1,\ldots,N_1. 
}
Notice that condition \eqref{eq:Christoffel_cond} holds for these measures. Specifically, we have
\bes{
\frac{1}{M_1} \sum^{M_1}_{i=1} \D \mu_i(\bm{y}) = {\cK(P_1,L^2(Z_0,\tau_0))(\bm{y})}\D\tau_0(\bm{y}),\quad \bm{y} \in D.
}
Hence, these sampling measures are suitable (in the sense of \eqref{L2_equiv}) for approximation over the initial domain estimate $\Omega_0$.

Now consider step $l = 2$ of GDAS. We have a new grid $Z_1$ (computed in stage (d) of step $l = 1$) and a new subspace $P_2$. To apply the ASGD methodology, we first restrict $\tau$ to $Z_1$. This gives the discrete probability measure
\bes{
\D \tau_1(\bm{y}) = \frac{\bbI_{Z_1}(\bm{y})}{\int_{Z_1} \D \tau(\bm{y})} \D \tau(\bm{y}),\quad \bm{y} \in D.
}
Notice that $Z_1$ is nonempty by construction and $\tau$ is supported on $Z$. Hence the denominator is nonvanishing, and therefore $\tau_1$ is well defined.
We now orthogonalize $P_2$ with respect to this measure, and write
\bes{
P_2 = \spn \{ \phi^{(2)}_1,\ldots,\phi^{(2)}_{N_2} \} \subset L^2(Z_1,\tau_1),
}
where $\{ \phi^{(2)}_i \}$ is an orthonormal basis for $P_2$, when considered as a subspace of $L^2(Z_1,\tau_1)$. Using this, we then define the new measures $\mu_{M_1+1},\ldots,\mu_{M_2}$ by
\bes{
\D\mu_{i}(\bm{y})=|\phi^{(2)}_{j}(\bm{y})|^2\D\tau_1(\bm{y}),\quad \bm{y}\in D,\ j = 1,\ldots,N_2,
}
where $i$ satisfies
\eas{
(j-1)(k_2 - k_1) + M_1  < i \leq j (k_2-k_1) + M_1,\qquad &j = 1,\ldots,N_1,
\\
(j-1) k_2 < i \leq j k_2,\qquad &j = N_{1}+1,\ldots,N_2.
}
Step $l \geq 2$ of GDAS proceeds in a similar manner. We have a new grid $Z_{l-1}$(computed in stage (d) of the previous step) and a new subspace $P_l$. We restrict $\tau$ to $Z_{l-1}$, giving a discrete probability measure
\be{
\D \tau_{l-1}(\mb{y}) 
= \frac{\bbI_{\newOmega_{l-1}}(\bm{y})}{\int_{\newOmega_{l-1}} \D \tau(\bm{y}) } \D \tau(\bm{y}),\ \bm{y}\in D, 
}
which is once again well defined. We now orthogonalize $P_l$ with respect to the measure $\tau_{l-1}$, and write
\bes{
P_l = \spn\{\phi_1^{(l)},\ldots,\phi_{N_l}^{(l)}\}\subset L^2(Z_{l-1},\tau_{l-1}),
}  
where $\{ \phi^{(l)}_i \}_{i=1}^{N_l}$ is the corresponding orthonormal basis for $P_l$, when considered as a subspace of $L^2(\newOmega_{l-1},\tau_{l-1})$.  Then we define the new measures $\mu_{M_{l-1}+1},\ldots,\mu_{M_{l}}$  
as
\bes{
\D\mu_i(\bm{y}) = |\phi_j^{(l)}(\bm{y})|^2\D\tau_{l-1}(\bm{y}),\ 
\bm{y}\in D,\ j = 1,\ldots,N_{l},
}
where $i$ satisfies
\be{
\label{i_j_relation}
\begin{split}
(j-1)(k_{l}-k_{l-1}) + M_{l-1}<i\leq j(k_l-k_{l-1}) + M_{l-1},\qquad & j=1,\ldots,N_{l-1},
\\
(j-1)k_l<i\leq j k_l,\qquad &  j=N_{l-1} + 1,\ldots,N_l.
\end{split}
}
 
\rem{ \label{rem:QR}
As with ASGD, in ASUD we use the orthonormal basis $\{ \phi^{(l)}_i \}^{N_{l}}_{i=1}$ for $P_l$ to construct the sampling measures. We once more compute this basis using QR decomposition. To do this, we first write the measure $\tau_{l-1}$ as
\bes{
\D\tau_{l-1}(\bm{y})=\sum_{i=1}^{K_{l-1}}\tau_i^{(l-1)}\delta(\bm{y}-\bm{z}_i^{(l-1)}),
}
where $\bm{z}_i^{(l-1)}$ is the $i$th point in the grid $Z_{l-1}=\{\bm{z}_i^{(l-1)}\}_{i=1}^{K_{l-1}}$ and  $\tau_i^{(l-1)}=\bbP(\bm{y}=\bm{z}_i^{(l-1)})$ for $\bm{y}\sim\tau_{l-1}$. Now, let $\{\psi_1,\ldots,\psi_{N_l}\}$ be a nonorthogonal basis for $P_l$ and
\bes{
\bm{B}_l = \left\lbrace 
\sqrt{\tau_{i}^{(l-1)}} \psi_j(\bm{z}_i^{(l-1)})
\right\rbrace_{i,j=1}^{K_{l-1},N_l} 
\in\bbC^{K_{l-1}\times N_l}.
}
Let $\bm{B}_l$ have a QR decomposition $\bm{B}_l=\bm{Q}_l\bm{R}_l$, where $\bm{Q}_l\in\bbC^{K_{l-1}\times N_l}$ and $\bm{R}_l\in\mathbb{C}^{N_l\times N_l}$, and write $\bm{Q}_l = \{q_{ij}^{(l)}\}_{i,j=1}^{K_{l-1},N_l}$. Then, as before,
\be{
\label{phi-psi-relation}
\phi_i^{(l)}(\bm{y}) = \sum_{j=1}^{i}(\bm{R}_l^{-\top})_{ij}\psi_j(\bm{y}),
\quad i=1,\ldots,N_l,
} 
is an  orthonormal basis for $P_l$ in $L^2(Z_{l-1},\tau_{l-1})$. Therefore, for $j = 1,\ldots,N_l$ and $i$ satisfying \R{i_j_relation}, we can write the sampling measure $\mu_i$ as
\eas{
\D\mu_i(\bm{y})=|\phi_j^{(l)}(\bm{y})|^2\D\tau(\bm{y})  
& = \sum_{k=1}^{K_{l-1}}\tau^{(l-1)}_k|\phi^{(l)}_j(\bm{z}_k^{(l-1)})|^2\delta(\bm{y}-\bm{z}_k^{(l-1)})\D\bm{y}
\\
& = \sum_{k=1}^{K_{l-1}}|q_{kj}^{(l)}|^2\delta(\bm{y}-\bm{z}^{(l-1)}_k)\D\bm{y}.
}
Hence $\mu_i$ is the discrete probability measure with $\bm{y}\sim\mu_i$ if $\bbP(\bm{y}=\bm{z}^{(l-1)}_k)=|q_{kj}^{(l)}|^2$ for $k=1,\ldots,K_{l-1}$. Sampling from $\mu_i$ is equivalent to drawing an integer $t$ randomly from $\{1,\ldots,K_{l-1}\}$ based on distribution $\{|q_{ij}^{(l-1)}|^2\}_{k=1}^{K_{l-1}}$ and then setting $\bm{y}=\bm{z}_t^{(l-1)}$.} 
 
With this remark in hand, we have now fully described the ASUD for GDAS. The resulting procedure is summarized in Method 2.   

\begin{tcolorbox}[title=Method 2.\ GDAS with Adaptive Sampling for Unknown Domains (ASUD)]
\noindent \textbf{Inputs:} Finite grid $Z$, probability measure $\tau$ over $Z$, subspace dimensions $0 = N_0 < N_1 < N_2 < \cdots$ and sampling ratios $1 \leq k_1 \leq k_2 \leq \cdots$ such that $M_l = k_l N_l$ for all $l$, function $\cQ$ as in \ref{def:Omega}. 
 
\noindent \textbf{Initialize:} Set $K_0=|Z|$, $\newOmega_0 = \{ \bm{z}^{(0)}_{i} \}^{K_0}_{i=1} = Z $, $I_0=\{1,\ldots,K_0\}$, $F_0 = 0$, $\cS_1 = \cR_1 = \emptyset$.
%\tcblower
\\
\textbf{For} $l=1,2,\ldots$ \textbf{do}\\
\begin{tcolorbox} 
\textbf{Stage (a)}
\tcblower
\textbf{(i)} Unless already defined, construct an approximation space $P_{l} = \spn \{ \psi_{j}\}^{N_l}_{j=1} \subset L^2(D,\tau)$ of dimension $N_{l}$ and set $M_{l}=k_{l}N_{l}$.\\
\textbf{(ii)} Construct the matrix $\bm{B}_{l}=\{\sqrt{\tau_i^{(l-1)}}\psi_j(\bm{z}^{(l-1)}_i)\}_{i,j=1}^{K_{l-1},N_{l}}$.  \\
\textbf{(iii)} Compute the QR decomposition $\bm{B}_{l}=\bm{Q}_{l}\bm{R}_{l}$ and write $\bm{Q}_{l}=\{q_{ij}^{(l)}\}_{i,j=1}^{K_{l-1},N_{l}}$.\\
\textbf{return} The discrete distributions $\{|q_{ij}^{(l)}|^2\}_{i=1}^{K_{l-1}}$ for $j = 1,\ldots,N_l$.
\end{tcolorbox}  

\begin{tcolorbox} 
\textbf{Stage (b)}
\tcblower
\textbf{if} $l>1$ \textbf{do}\\
\textbf{for} $t = 1,\ldots,N_{l-1}$ \textbf{do}\\
\textbf{for} $s=1,\ldots,k_l-k_{l-1}$ \textbf{do}\\
Draw an integer $v$ randomly and independently from $I_{l-1}$ according to the discrete distribution $\{ |q_{it}^{(l-1)}|^{2} \}^{K_{l-1}}_{i=1}$. If $\cQ(f(\bm{z}_v))=0$ then add $\bm{z}_v$ to $\cR_l$. Repeat until $\cQ(f(\bm{z}_v)) = 1$ and then add $\bm{z}_v$ to $\cS_l$.
\\
\textbf{end for} \\
\textbf{end for} \\
\textbf{end if}  \\
\textbf{for} $t = N_{l-1}+1,\ldots,N_l$ \textbf{do}\\
\textbf{for} $s = 1,\ldots,k_l$ \textbf{do}\\
Draw an integer $v$ randomly and independently from $I_{l-1}$ according to $\{ |q_{it}^{(l-1)}|^{2} \}^{K_{l-1}}_{i=1}$. If $\cQ(f(\bm{z}_v))=0$ then add $\bm{z}_v$ to $\cR_l$. Repeat until $\cQ(f(\bm{z}_v)) = 1$ and then add $\bm{z}_v$ to $\cS_l$.
\\
\textbf{end for}
\\ 
\textbf{end for}
\\
\textbf{return} Sample points $\cS_l = \{ \bm{y}_i \}^{M_l}_{i=1}$ and rejected points $\cR_l = \{ \bm{u}_i \}^{T_l}_{i=1}$, total number of function evaluations $F_l = F_{l-1} + G_l$, where $G_l$ is the number of evaluations of $f$ used in this stage, and $T_l = F_l - M_l$.
\end{tcolorbox}
\textbf{Stage (c)} Compute an approximation $\tilde{f}_{l}$ to $f$ from the data $\{f(\bm{y}_i)\}^{M_l}_{i=1}$.
\\
\textbf{Stage (d)}  Compute $\newOmega_l = \{ \bm{z} \in Z : \cQ(\tilde{f}_l(\bm{z})) = 1 \} \cup \cS_l \setminus \cR_l$. Define a new set of indices $I_l=\{j_1,\ldots,j_{K_{l}}\}$ so that $\newOmega_l = \{ \bm{z}_{j_i} \}^{K_l}_{i=1} : =   \{ \bm{z}^{(l)}_{i} \}^{K_l}_{i=1}$, where $K_l = |Z_l|$.
\\
Update $\cS_{l+1} = \cS_l$ and $\cR_{l+1} = \cR_l$
\\
\textbf{end for}
\\
\textbf{Output:} The sequences of function approximations $\{\tilde{f}_l\}_{l\geq 1}$ and domain approximations $\{Z_l\}_{l\geq 1}$.
\end{tcolorbox}

\section{Approximation methods for $\tilde{f}$} 
\label{sec:other_decoders}

Note that ASUD allows for arbitrary methods of approximation in stage (c). We now discuss the two approaches for doing this that we consider later in our numerical experiments. Both approaches assume that the subspaces $P_l$ are defined \textit{a priori} as a sequence of nested finite-dimensional subspaces $P_1\subseteq P_2\subseteq \ldots\subset L^2(D,\tau)$. Note that other approaches are possible within our framework, including those that adaptively generate the $P_l$. See \S \ref{sec:conclusions} for some further discussion on this point.

\subsection{(Weighted) least-squares approximation}
\label{sec:least_sq}

Let $N_l = \dim(P_l)$ and suppose that $M_l \geq N_l$. Let $\bm{y}_1,\ldots,\bm{y}_{M_l}$ be the sample points generated up to and including stage (b) of step $l$ of ASUD. Let $\tau_l$, $\mu_l$ and $\newOmega_{l-1}$ be as in \S \ref{sec:ASUD}, and let
\bes{
w_l= \left ( {\cK(P_l,L^2(Z_{l-1},\tau_{l-1}))} \right )^{-1},
} 
where ${\cK(P_l,L^2(Z_{l-1},\tau_{l-1}))}$ is the reciprocal Christoffel function of $P_l$ as a subspace of $L^2(Z_{l-1},\tau_{l-1})$ (see \R{Christoffel_fun}). Then, similar to in Remark \ref{r:samplecomplexitystable_LS}, we define the (weighted) least-squares approximation of $f$ as 
\be{\label{formulation:wls}
\tilde{f}_l \in 
\argmin{p \in P_l} 
\left\lbrace 
\frac{1}{M_l}\sum_{i=1}^{M_l}
w_l(\bm{y}_i) |f(\bm{y}_i)-p(\bm{y}_i)|^2
\right\rbrace.
} 
Write $\tilde{f}_l=\sum_{i=1}^{N_l}c_i^{(l)}\phi_i^{(l)}$ in terms of the orthonormal basis $\{ \phi^{(l)}_{j} \}^{N_l}_{j=1}$ for the subspace $P_l \subset L^2(Z_{l-1},\tau_{l-1})$. Then the coefficients $c^{(l)}_i$ of $\tilde{f}_l$ are a solution of the algebraic least-squares problem
\be{
\label{alg_LS_soln}
\bm{c}^{(l)} = (c_i^{(l)})_{i=1}^{N_l}\in \argmin{\bm{x}\in\mathbb{C}^{N_l}}\nm{\bm{A}^{(l)} \bm{x}-\bm{b}^{(l)} }_2,
}
where 
\bes{
\bm{A}^{(l)}=
\left\lbrace
\sqrt{\frac{w_l(\bm{y}_i)}{M_l}}
\phi_j^{(l)}(\bm{y}_i)
\right\rbrace_{i=1,j=1}^{M_l,N_l}
\in \mathbb{C}^{M_l\times N_l},
\quad
\bm{b}^{(l)}=
\left\lbrace 
\sqrt{\frac{w_l(\bm{y}_i)}{M_l}} f(\bm{y}_i)
\right\rbrace_{i= 1}^{M_l}
\in \mathbb{C}^{M_l}.
}
Now, as in stage (a)(iii) of ASUD (Method 2), let $\bm{Q}_l$ be the matrix arising from the QR factorization of $\bm{B}_l$, with entries $\{ q^{(l)}_{ij} \}^{K_{l-1},N_l}_{i,j=1}$. Since the points $\bm{y}_1,\ldots,\bm{y}_{M_l}$ belong to $Z_{l-1}$, we can write $\bm{y}_{i} = \bm{z}^{(l)}_{j_i}$, where $j_i \in \{1,\ldots,K_{l-1} \}$ for $i = 1,\ldots,M_l$. Hence, we can rewrite $\bm{A}^{(l)}$ and $\bm{b}^{(l)}$ as
\bes{
\bm{A}^{(l)}=\left\lbrace \frac{q_{j_i k}^{(l)}}{\sqrt{\frac{M_l}{N_l}\sum_{t=1}^{N_l}|q_{j_i t}^{(l)}|^2}}\right\rbrace_{i,k=1}^{M_l,N_l},\quad 
\bm{b}^{(l)}=\left\lbrace \frac{f(\bm{y}_{i})}{\sqrt{\frac{M_l K_{l-1}}{N_l}\sum_{t=1}^{N_l}|q_{j_i t}^{(l)}|^2}}\right\rbrace_{i=1}^{M_l},
}
in terms of the matrix $\bm{Q}_{l}$. In particular, $\bm{A}^{(l)}$ consists of rows of $\bm{Q}_l$ corresponding to the indices $j_1,\ldots,j_{M_l}$, scaled by values $\sqrt{w_l(\bm{y}_i) / M_l}$.
 
Let $\bm{c}^{(l)}$ be a solution \R{alg_LS_soln} of the algebraic least-squares problem. In order to perform stage (d) of ASUD, we need to evaluate the corresponding function $\tilde{f}_l$ over the grid $Z_0$. This can be done straightforwardly using the matrix $\bm{R}_{l}$ arising from the QR factorization of $\bm{B}_{l}$, as well as the matrix 
\bes{
\bm{C}_l = \left \{ \psi_j ( \bm{z}^{(0)}_i ) \right \}^{K_0,N_l}_{i,j = 1}. 
}
Specifically, the vector
\bes{
\tilde{\bm{f}}_l = \bm{C}_{l} (\bm{R}_l)^{-1}\bm{c}^{(l)},
}
contains the values of $\tilde{f}_l$ over $Z_0$. Indeed, the action of $(\bm{R}_{l})^{-1}$ on $\bm{c}^{(l)}$ yields the coefficients of $\tilde{f}_l$ in the nonorthogonal basis $\{ \psi_i \}^{N_l}_{i=1}$ for $P_l$, and the subsequent action of $\bm{C}_l$ evaluates an expansion in this basis (given its coefficients) on the grid $Z_0$. 
Indeed, consider $\bm{z}^{(0)}_t\in Z_0$ and \R{phi-psi-relation}. Then
\eas{
\tilde{f}_l(\bm{z}^{(0)}_t) = \sum_{i=1}^{N_l}c_i^{(l)}\phi_i^{(l)}(\bm{z}^{(0)}_t) = \sum_{i=1}^{N_l}c_i^{(l)}\sum_{j=1}^{i}(\bm{R}_l)^{-\top}_{ij}\psi_j(\bm{z}^{(0)}_t) & = \sum_{i=1}^{N_l}c_i^{(l)}\sum_{j=1}^{i}(\bm{C}_l)_{tj}(\bm{R}_l)^{-\top}_{ij},
}
which is precisely $(\bm{C}_l(\bm{R}_l)^{-1}\bm{c}^{(l)})_t$, as required.
Note that computing $\tilde{\bm{f}}_l$ is not only necessary for stage (d) of ASUD, it also allows us to represent $\tilde{f}_l$ over the domain approximation $Z_{l}$ and, in simulations, the true (discrete) domain $\Omega \cap Z$. Hence, instead of the function $\tilde{f}_l$, we consider the vector $\tilde{\bm{f}}_l$ as the output of the weighted least-squares approximation. This approximation is summarized below:

\begin{tcolorbox}[floatplacement=th!,float,title= (Weighted) least squares]
\textbf{Inputs:} $\bm{C}_l$, QR decomposition $\bm{Q}_l \bm{R}_l$ of $\bm{B}_l$, and indices $i_1,\ldots,i_{M_l}$. 
\\
\textbf{(i)} Define $\bm{A}^{(l)}\in\bbC^{M_l\times N_l}$ and $\bm{b}^{(l)}\in\bbC^{M_l}$ as 
\bes{
\bm{A}^{(l)}=\left\lbrace \frac{q_{j_{i}k}^{(l)}}{\sqrt{\frac{M_l}{N_l}\sum_{t=1}^{N_l}|q_{j_{i}t}^{(l)}|^2}}\right\rbrace_{i,k=1}^{M_l,N_l}
\text{ and }
\bm{b}^{(l)}=\left\lbrace \frac{f(\bm{y}_{i})}{\sqrt{\frac{M_l K_{l-1}}{N_l}\sum_{t=1}^{N_l}|q_{j_{i}t}^{(l)}|^2}}\right\rbrace_{i=1}^{M_l}.
}
\textbf{(ii)} Compute $\bm{c}^{(l)}=\argmin{\bm{x}\in\bbC^{N_l}}\nm{\bm{A}^{(l)}\bm{x}-\bm{b}^{(l)}}_{2}$ and $\tilde{\bm{f}}_l=\bm{C}_l(\bm{R}_l)^{-1}\bm{c}_l$.
\\ 
\textbf{return} Approximation $\tilde{\bm{f}}_l$ to $f$.
 
\end{tcolorbox}

\subsection{Augmented (weighted) least-squares approximation}
\label{sec:aug_least_sq}

A limitation of the approximation described previously is that it makes no use of the rejected points $\cR_l$. 
There are a number of ways that one might strive to incorporate these points into an approximation scheme. We now describe a modification of the previous procedure that can lead to significantly better domain learning for certain functions.

As before, we consider a sequence of nested finite-dimensional subspaces $P_1\subseteq P_2\subseteq \ldots\subset L^2(D,\tau)$ with $N_l = \dim(P_l)$. We now make several modifications to Method 2. First, we modify stage (b) to only reject a point $\bm{z}_v$ if $f(\bm{z}_v) = + \infty$. In other words, we accept points both inside $\Omega$ and outside $\Omega$. Then, in Stage (d) we modify the domain update as follows:
\eas{
Z_l = & \left ( \{ \bm{z} \in Z : \cQ(\tilde{f}_l(\bm{z})) = 1 \} \cup \{ \bm{y}_i : \cQ(f(\bm{y}_i)) = 1,\ i = 1,\ldots,M_l \} \right ) 
\\
& \Big\backslash \left ( \{ \bm{y}_i :  \cQ(f(\bm{y}_i)) = 0,\ i = 1,\ldots,M_l \} \cup \{ \bm{u}_i \}^{T_l}_{i=1} \right ).
}
In other word, our domain estimate includes all points in $Z$ for which  $\cQ(\tilde{f}_l(\cdot)) = 1$ and all the sample points $\bm{y}_i$ which belong to $\Omega$, and excludes those sample points that do not belong to $\Omega$ and any rejected points $\bm{u}_i$ (i.e.\ those for which $f = + \infty$).

Given these modifications, we define the resulting approximation exactly as before, via weighted least squares. Since it uses points both inside and outside of $\Omega$, we term this procedure \textit{augmented (weighted) least-squares approximation}. 
 
\section{Numerical experiments} 
\label{sec:num_exp}

In this section, we present numerical experiments demonstrating the performance of ASUD with both unaugmented and augmented weighted least-squares approximation.
 
\subsection{Experimental setup}
\label{sec:num_setup}
In this section, we describe our experimental setup. Throughout, we consider the domain $D = [-1,1]^d$. We let $Z_0 = Z$ be a uniform grid of size $K=30000$, drawn uniformly and randomly from $D$, and let $\tau$ be the uniform measure over $Z$. The grid $Z$ is generated once before all the subsequent computations. Throughout, we consider the approximation of smooth functions using polynomials. To do so, we choose
\bes{
P_l = P^{\mathrm{HC}}_{n_l} = \spn \left\{ \bm{y} \mapsto \bm{y}^{\bm{n}} : \bm{n} \in \Lambda^{\mathrm{HC}}_{n_l} \right \},
}
where $\Lambda^{\mathrm{HC}}_{n}$ is the \textit{hyperbolic cross} index set of index $n$:
\bes{
\Lambda^{\mathrm{HC}}_n = \left \{ \bm{n} = (n_1,\ldots,n_d) \in \bbN^d_0 : \prod^{d}_{k=1} (n_k+1)\leq n+1 \right \}.
}
The initial basis $\{\psi_1,\ldots,\psi_{N_l}\}$ for $P_l$ is constructed by taking the restrictions to $Z_{l-1}$ of the orthonormal Legendre polynomials on $[-1,1]^d$ with indices in $\Lambda$.  We also define sampling rates $0=M_0<M_1<\ldots<M_r=M_{\max}$ as $M_l=k_l N_l$, where $N_l=\dim(P_l)$, $1\leq N_1<\ldots<N_r=N_{\max}\leq 1000$ and $k_l$ is the closest integer to $\log(N_l)$. 
 
\begin{figure}[t]
\begin{center} 
\begin{tabular}{ccc} 
\includegraphics[scale=0.18]{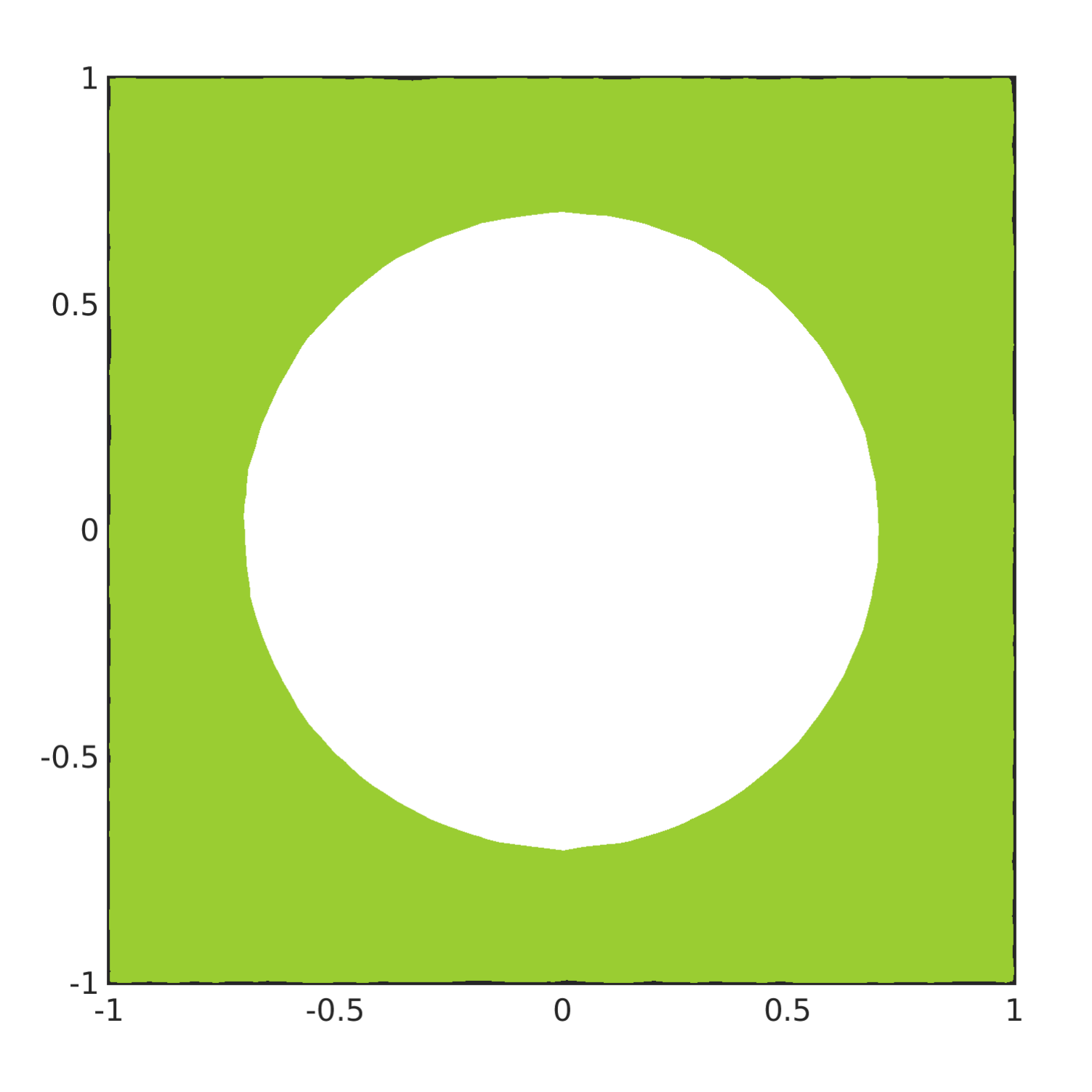} & 
\includegraphics[scale=0.18]{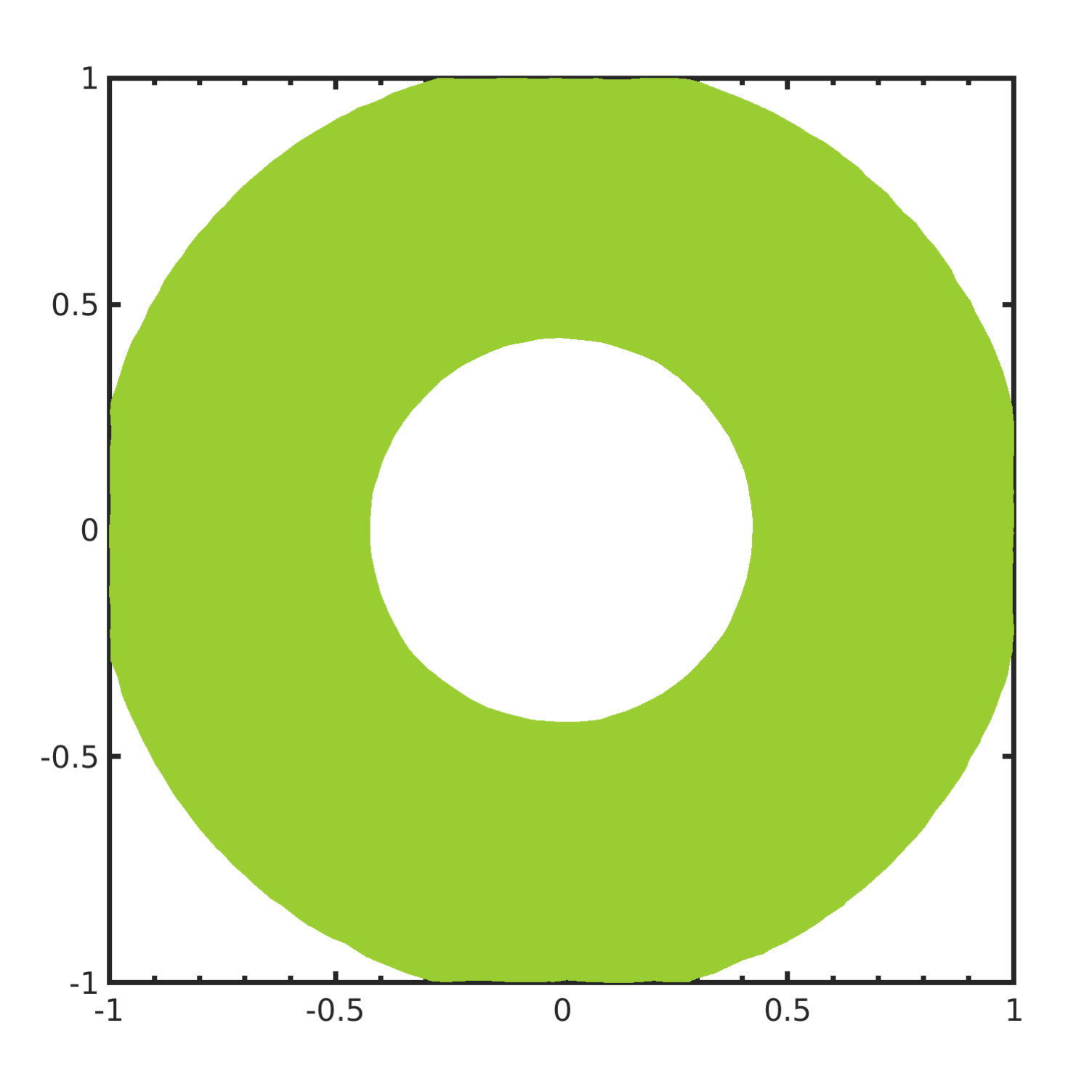} &  
\includegraphics[scale=0.18]{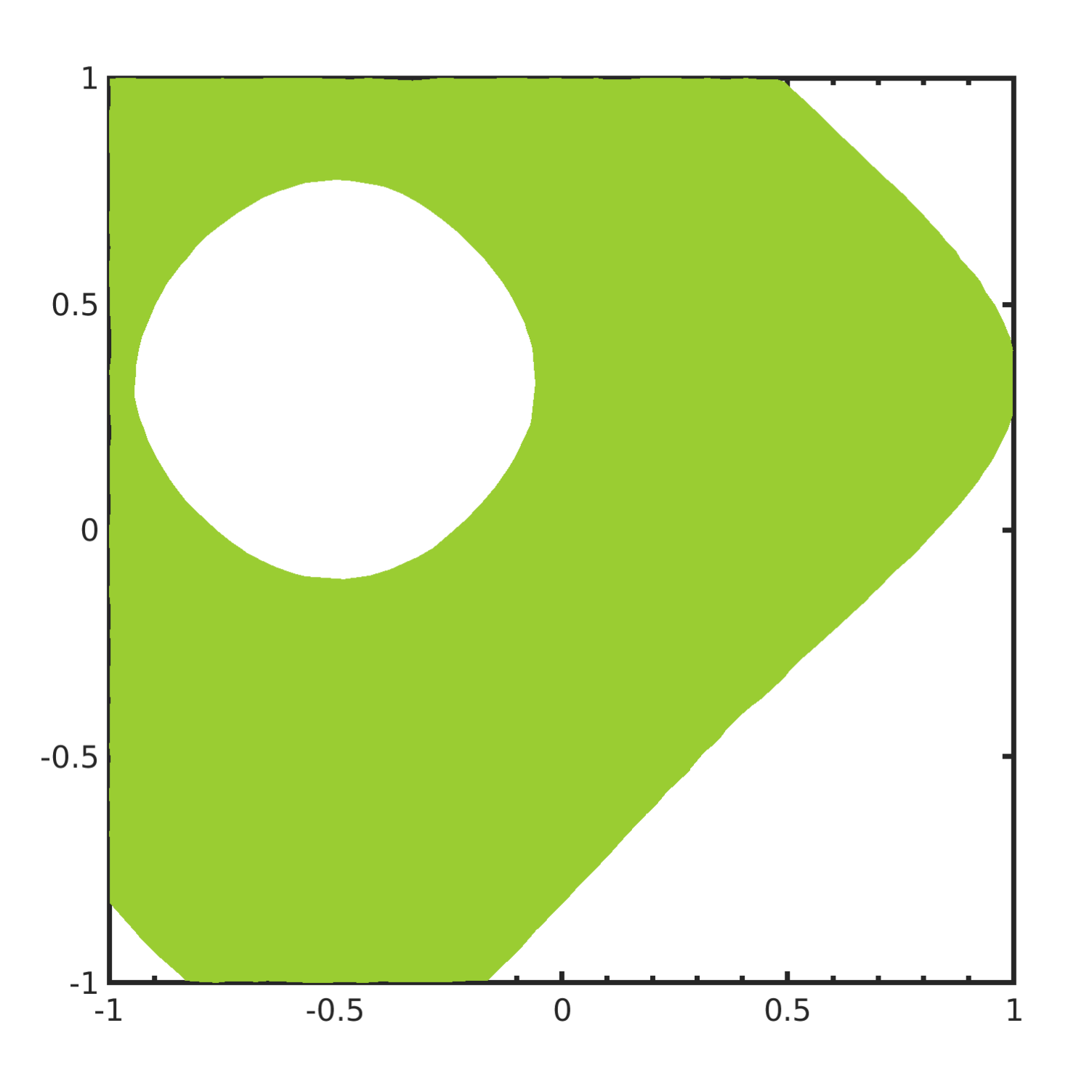} 
\end{tabular} 
\end{center}
\caption{
Domains $\Omega_1$, $\Omega_2$, $\Omega_4$, (left to right) for $d=2$. Note that $\Omega_2=\Omega_3$ when $d=2$.  
}
\label{fig:Omega_123}
\end{figure}

We consider the following test functions:  
\eas{ 
f_1(\bm{y})& =\left(\left(\frac{10}{7}\right)^2-\frac{1}{y_1^2+y_2^2}\right)\exp\left(-\sum_{i=1}^d y_i/2d\right), 
\\
f_2(\bm{y}) &= \log\left(8\sum_{i=1}^d y_i^2\right)-2\left(\sum_{i=1}^dy_i^2\right),   
\\
f_3(\bm{y}) &= \left(1-\frac{(d-2)}{100}(d^2-10d+29)\right)\log\left(\frac{16}{d}\sum_{i=1}^d y_i^2\right)-\frac{4}{d}\left(\sum_{i=1}^dy_i^2\right),
\\
f_4(\bm{y}) &= \prod_{i=1}^d\frac{d/4}{d/4 + (y_i + (-1)^{i+1}/(i+1))^2}.  
} 
We define the corresponding domain $\Omega = \Omega_i$ for the function $f_i$ as
\eas{
\Omega_i &=\{\bm{z}\in D: 0 \leq f_i(\bm{y}) < \infty\},\ i = 1,2,3,
\quad
\Omega_4 =  \{\bm{z}\in D: 0.18\leq f_4(\bm{z})\leq 0.72\}. 
} 
Examples of these domains in $d = 2$ dimensions are shown in Fig.\ \ref{fig:Omega_123}. Note that $f_2=f_3$ when $d=2$. All four functions are smooth (analytic) within their respective domains. Notice, however, that $f_2$ and $f_3$ have singularities at $\bm{y}= \bm{0}\in D$ and $f_1$ has a singularity at any $\bm{y} \in D$ for which $y_1 = y_2 = 0$. The reason for considering functions with singularities at certain points is that it allows us to model cases where the `black box' evaluating $f$ returns an exit flag (NaN, Inf), which we think of as $+\infty$.

To measure the error between $f$ and its approximation $\tilde{f}$, we consider the relative approximation error 
\be{\label{eq:error_l}
E_l(f)=\frac{\nmu{f-\tilde{f}_l}_{L^2(\Omega,\tilde{\tau})}}{\nm{f}_{L^2(\Omega,\tilde{\tau})}},\quad
\mathrm{d}\tilde{\tau}(\bm{y}):=\frac{\bbI_{\Omega}(\bm{y})}{\int_\Omega d\tau(\bm{y})}\mathrm{d}\tau(\bm{y}),
}
where $\tilde{\tau}$ is the restriction of $\tau$ to $\Omega$.
In order to measure the approximation of the true domain, we define the mismatch volume between $\Omega$ and $Z_l$ as 
\be{
V_l(f)
= \frac{|(Z_{\Omega} \setminus Z_{l-1})\cup(Z_{l-1}\setminus Z_{\Omega})|}{|Z_{\Omega}|},
}
where $Z_{\Omega} = \Omega \cap Z$ is the discrete representation of $\Omega$.
As we mentioned in Section \ref{sec:setup_discussion}, we also compute the \textit{rejection rate} as 
\be{
R_l(f) = \frac{F_l-M_l}{F_l},
} 
where $F_l$ is the number of function evaluations of ASUD at the $l$th step. Since our methods involve random sampling, we perform multiple trials and then average the corresponding quantities $E_l(f)$, $V_l(f)$ and $R_l(f)$. 
Throughout, we present the mean values of these quantities averaged over $50$ trials wtih respect to the sample points. 

We consider ASUD with either the least-squares approximation of Section \ref{sec:least_sq}, which we term ASUD-LS, or the augmented least-squares approximation of Section \ref{sec:aug_least_sq}, which we term ASUD-ALS. We also compare these schemes against two other methods. The first is the ASGD method of \cite{AdcockCardenas2020}, which was discussed earlier in Section \ref{alg:asgd}. We term this ASGD-LS. Note that this method assumes the domain $\Omega$ is known. 
Thus, it provides a benchmark against which to compare the performance of ASUD.  Since the number of function evaluations $F_l = M_l$ in this case, we do not report the rejection rate $R_l(f)$ for this method. At the other end of the spectrum, we also consider (unweighted) least-squares with Monte Carlo sampling, which we term MC-LS. As discussed in Section \ref{sec:intro}, this method is often used in practice. Yet it can be wasteful, since it does not adapt the sampling to a domain estimate. Later, we will see that this method often has a high rejection rate.

\subsection{Experimental results}

\begin{figure}[!h]
\begin{center}
\begin{tabular}{ccc}  
\includegraphics[scale=0.35]{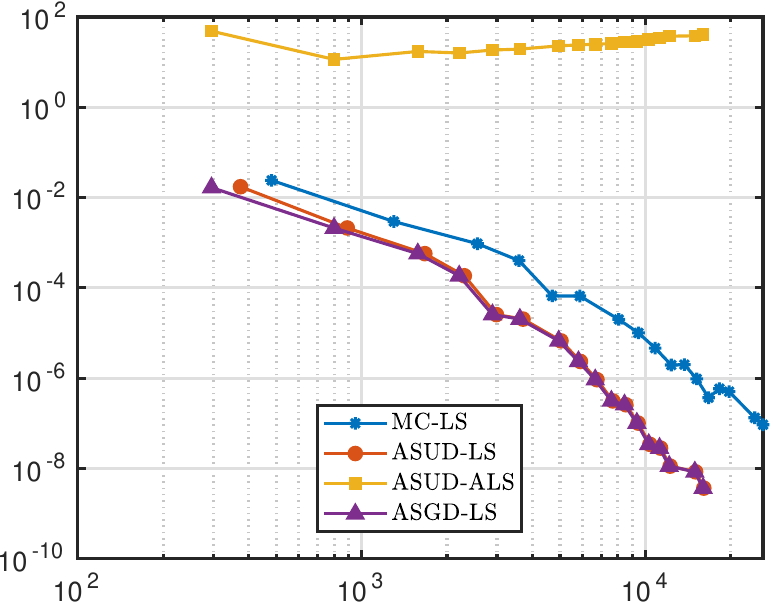} &  
\includegraphics[scale=0.35]{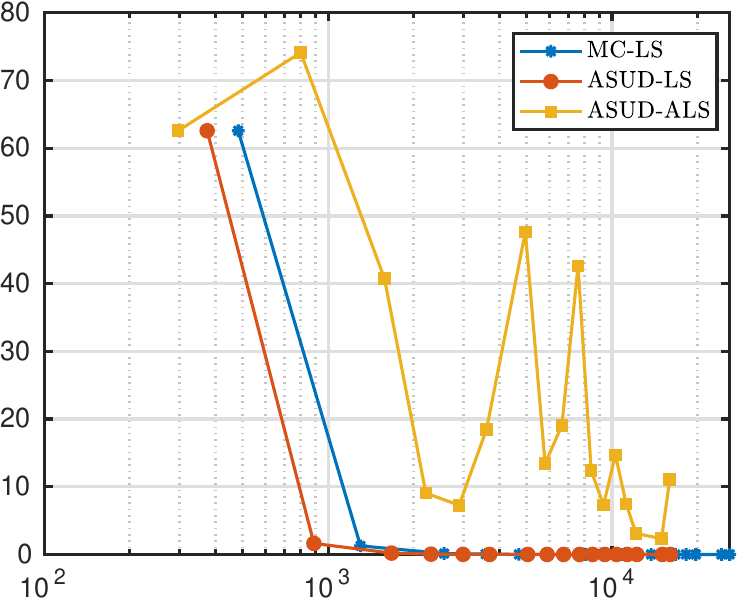} &  
\includegraphics[scale=0.35]{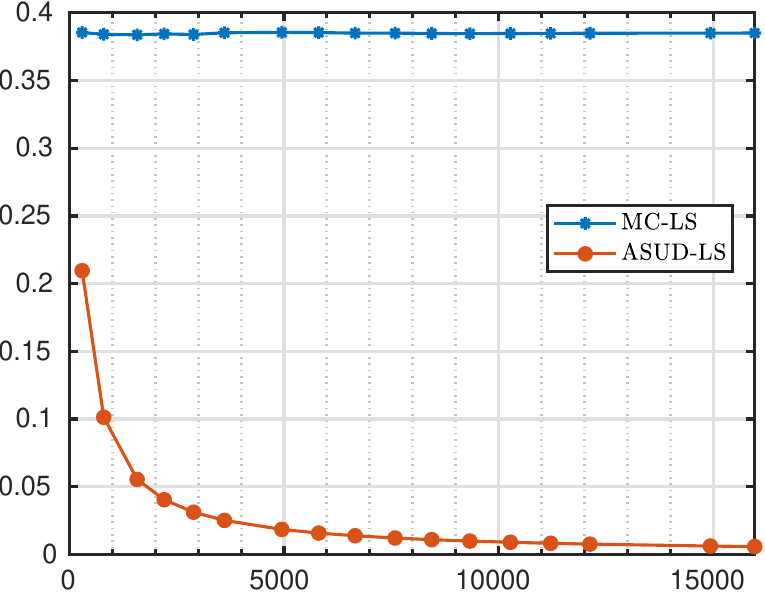} \\   
\includegraphics[scale=0.35]{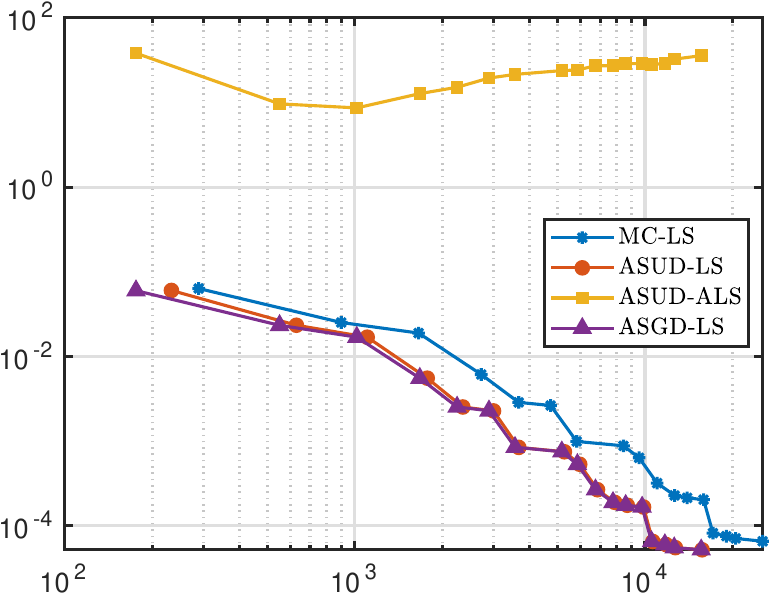} &  
\includegraphics[scale=0.35]{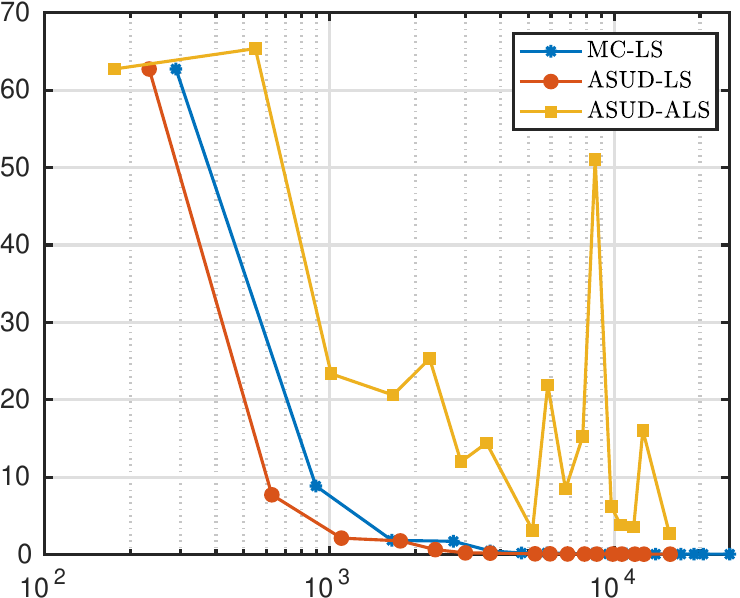} & 
\includegraphics[scale=0.35]{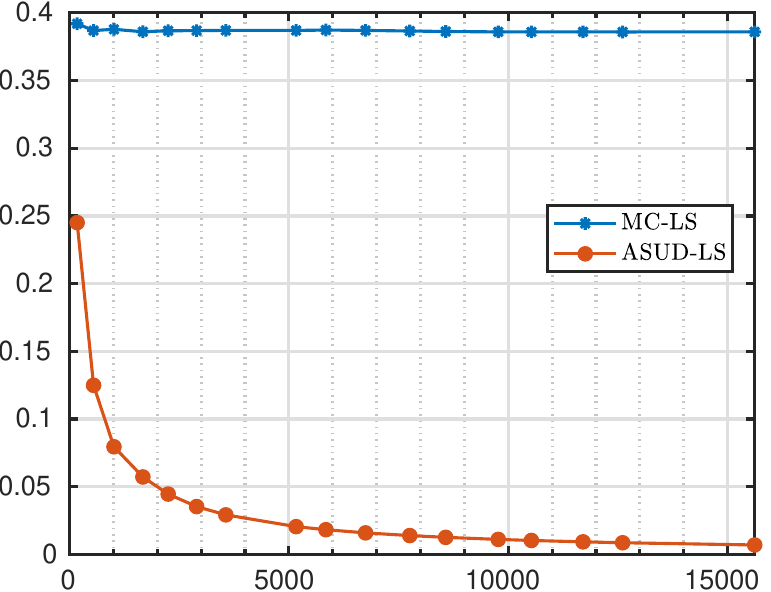} \\    
\includegraphics[scale=0.35]{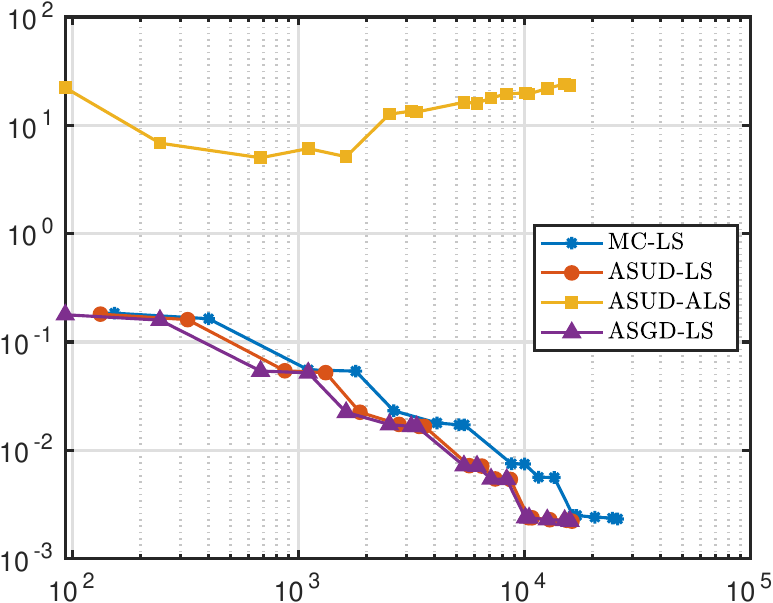} &  
\includegraphics[scale=0.35]{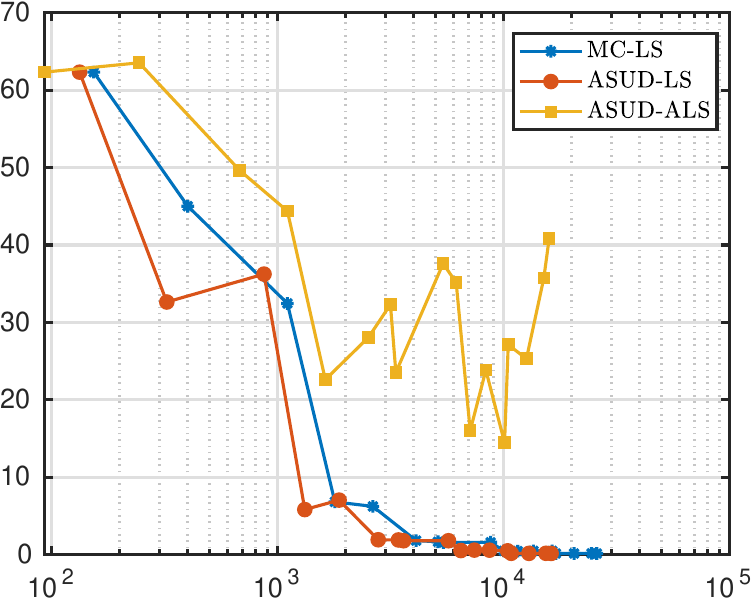} &  
\includegraphics[scale=0.35]{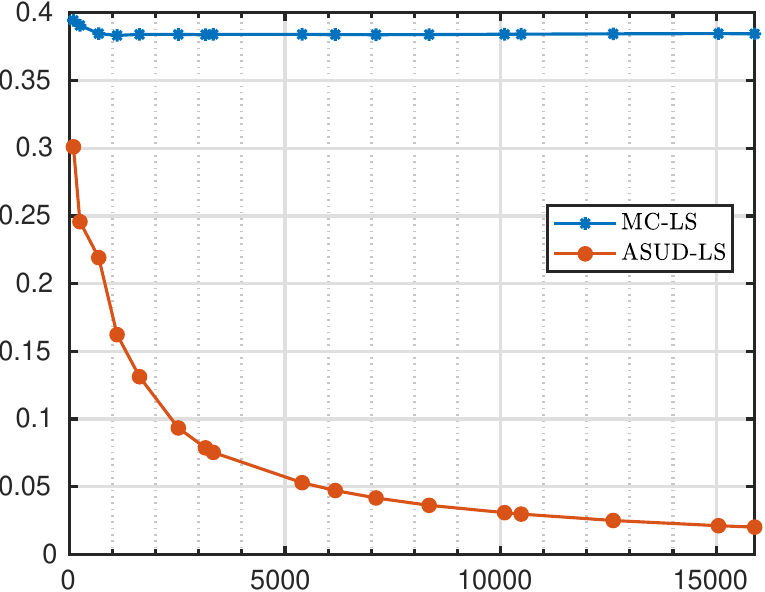} \\  
\includegraphics[scale=0.35]{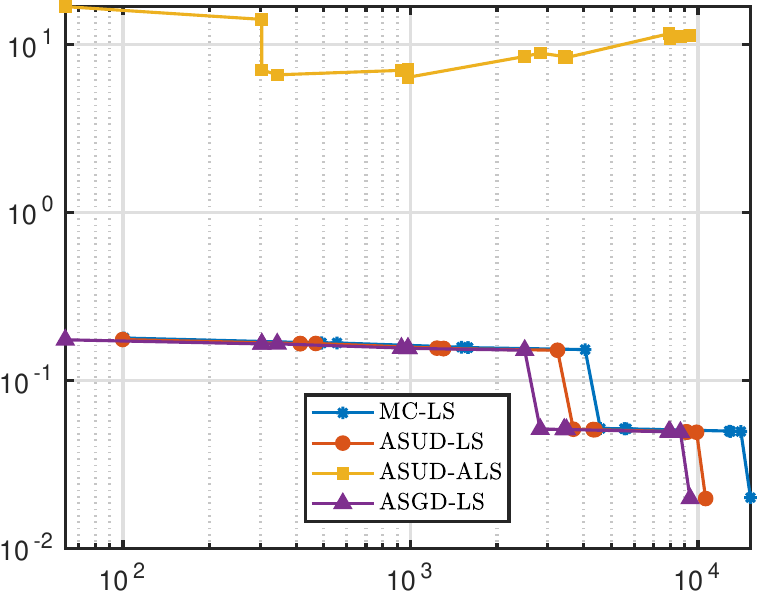} &  
\includegraphics[scale=0.35]{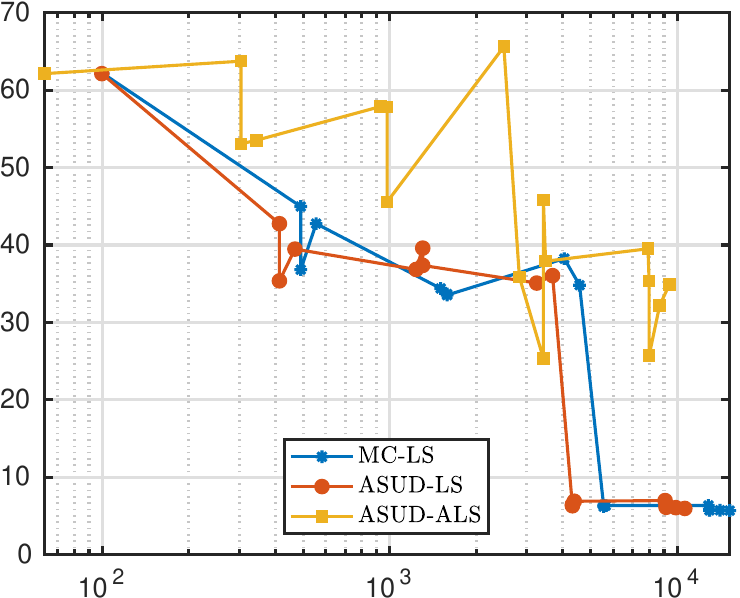} &  
\includegraphics[scale=0.35]{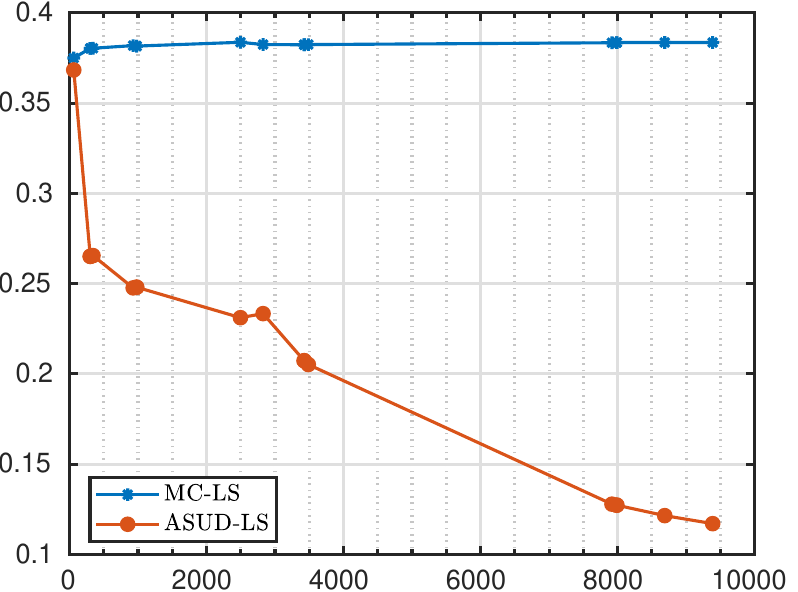} \\  
\includegraphics[scale=0.35]{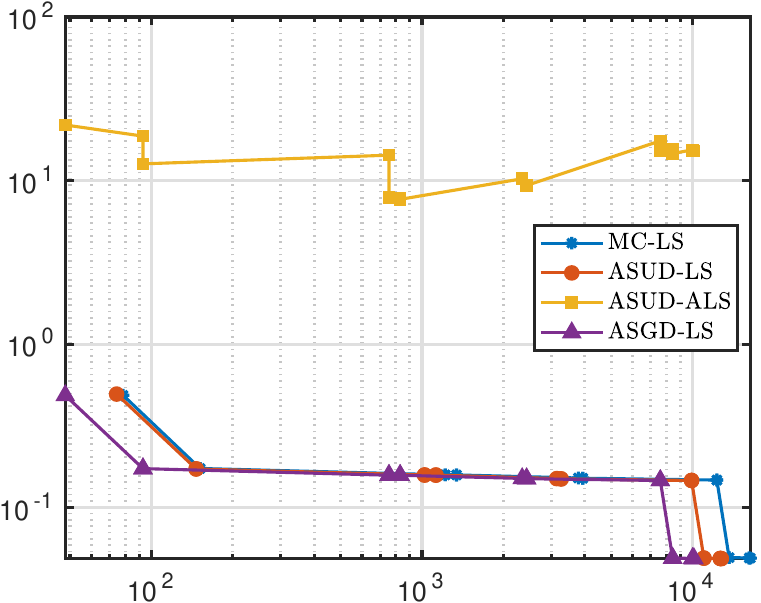} &  
\includegraphics[scale=0.35]{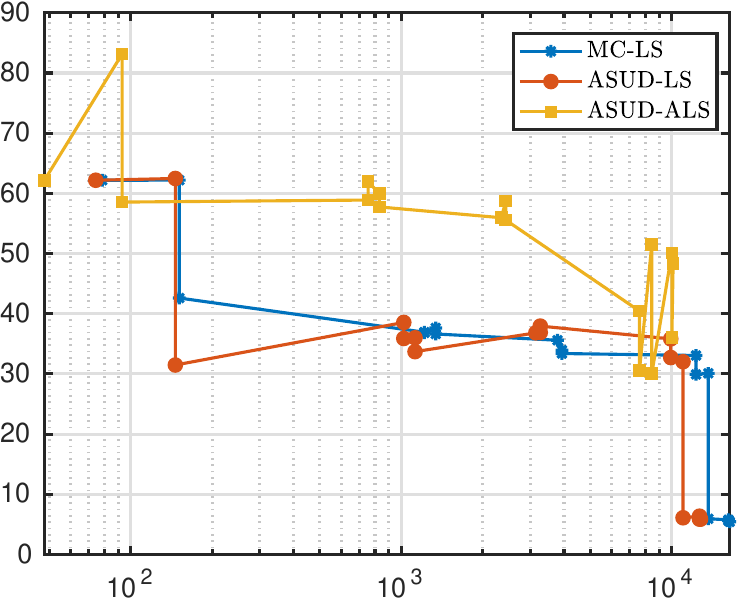} &  
\includegraphics[scale=0.35]{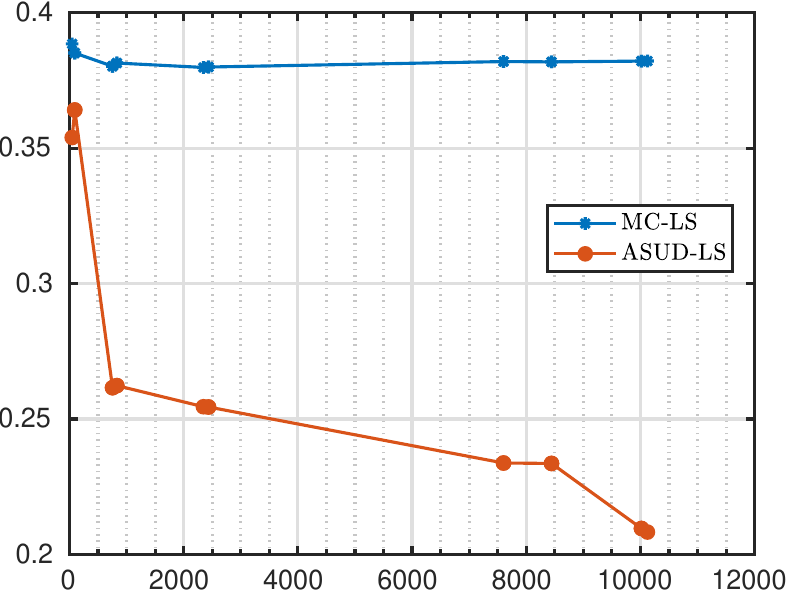} 
\end{tabular}
\end{center}
\vspace*{-1mm}
\caption{Approximation of the function $f = f_1$ and domain $\Omega = \Omega_1$ in $d = 2,3,5,10,15$ dimensions (top to bottom). Left: the relative error $E_l(f)$ versus the number of function evaluations $F_l$. Middle: the mismatch volume $V_l(f)$ versus $F_l$. Right: the rejection rate $R_l(f)$ versus $M_l$.
}  
\label{fig:F3_highdim}
\end{figure}

In Fig.\ \ref{fig:F3_highdim}, we show results for the function $f = f_1$. We see a clear benefit of ASUD-LS over MC-LS in lower dimensions with respect number of function evaluations needed to obtain a certain accuracy. In fact, the ASUD-LS approximation error is quite close to that of ASGD-LS, despite it assuming no \textit{a priori} knowledge of the domain. As shown in the right column of this figure, the rejection rate $R_l(f)$ for MC-LS is around $40\%$, since $\mathrm{Vol}(\Omega) / \mathrm{Vol}(D) \approx 0.6$ in all dimensions. Conversely, the rejection rate is decreasing for ASUD-LS, and close to zero in dimensions $d = 2,3,5$ for large enough $l$. Even in higher dimensions, however, the rejection rate is significantly smaller than for MC-LS. This translates to it needing fewer function evaluations to achieve a certain error. Indeed, when $d = 15$ we achieve the minimum error using roughly $2\times 10^3$ function evaluations versus $3\times 10^3$ function evaluations for MC-LS.
  
 It is notable that ASUD-ALS gives a significantly worse approximation than the other methods. We discuss the reasons for this later. In particular, the mismatch volume $V_l(f)$ for this method is never below $20\%$, while the other methods can achieve close to zero domain learning errors in the lower-dimensional cases. On the other hand, when $d = 15$ ASUD-LS achieves a similar domain learning error for large enough numbers of function evaluation, despite yielding a much worse approximation error.

\begin{figure}[!h]
\begin{center}
\begin{tabular}{ccc} 
\includegraphics[scale=0.35]{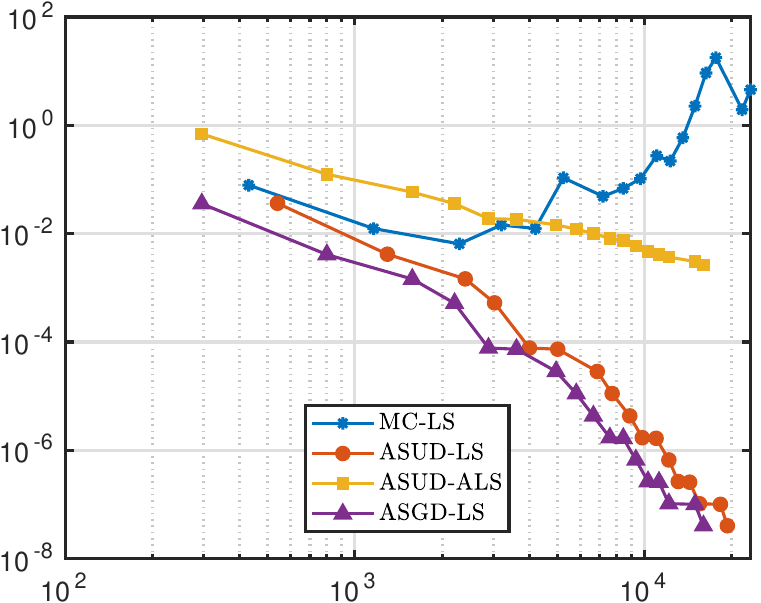} &  
\includegraphics[scale=0.35]{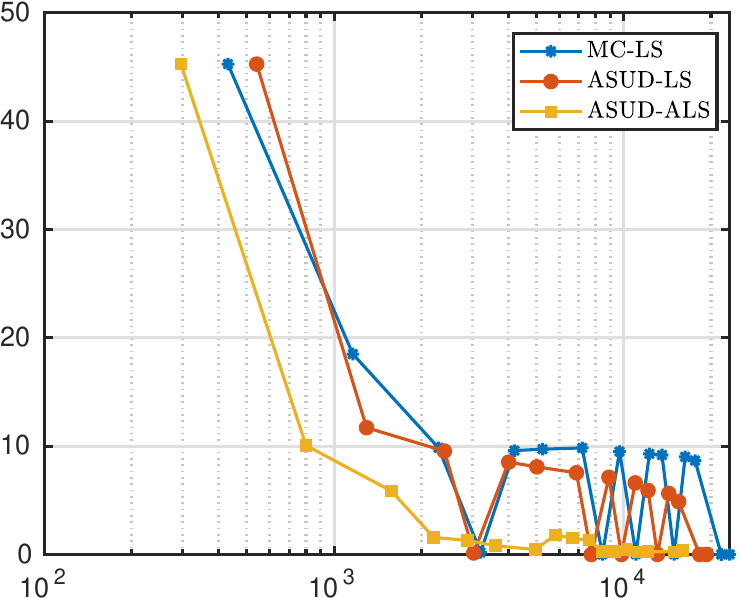} &  
\includegraphics[scale=0.35]{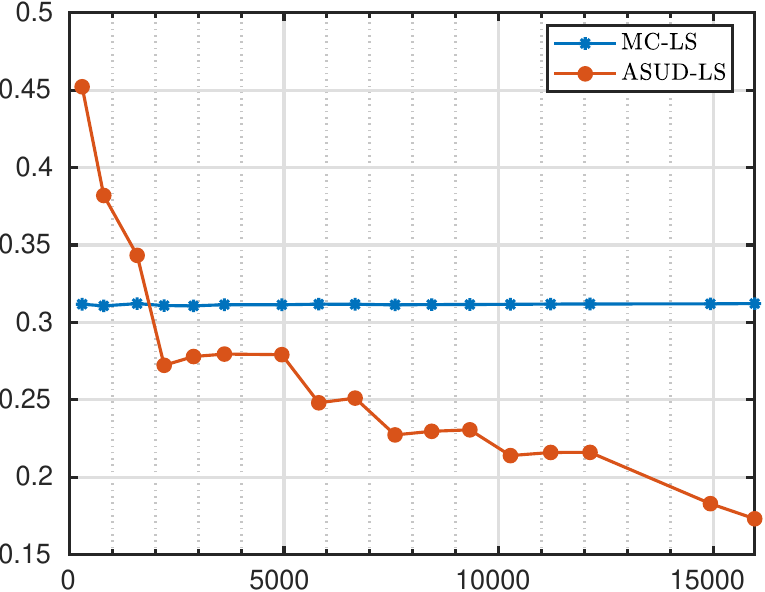} \\ 
\includegraphics[scale=0.35]{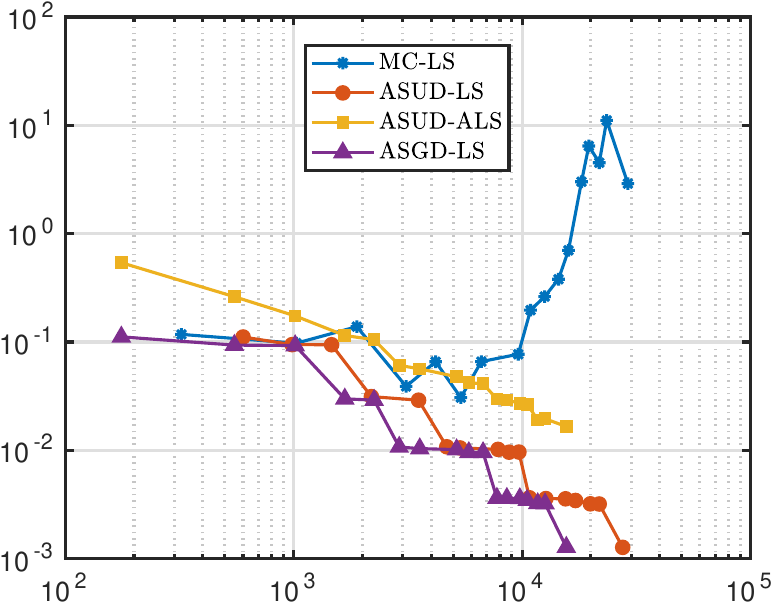} &  
\includegraphics[scale=0.35]{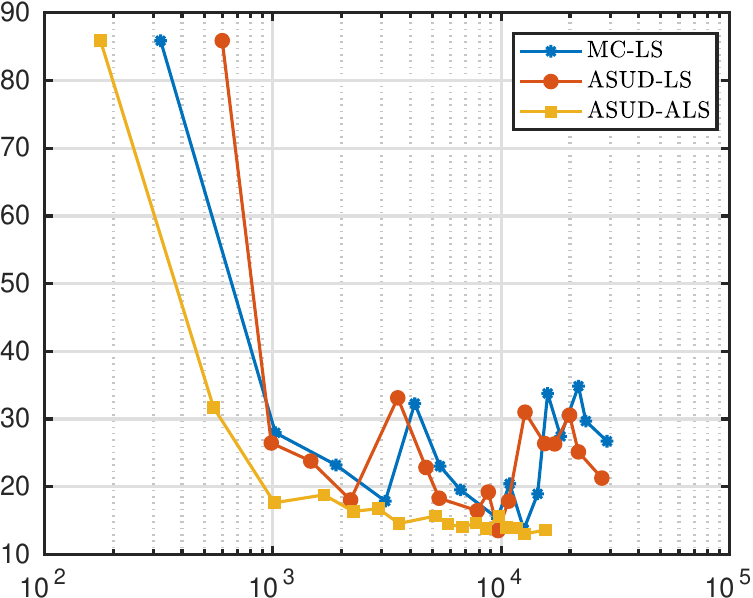} &  
\includegraphics[scale=0.35]{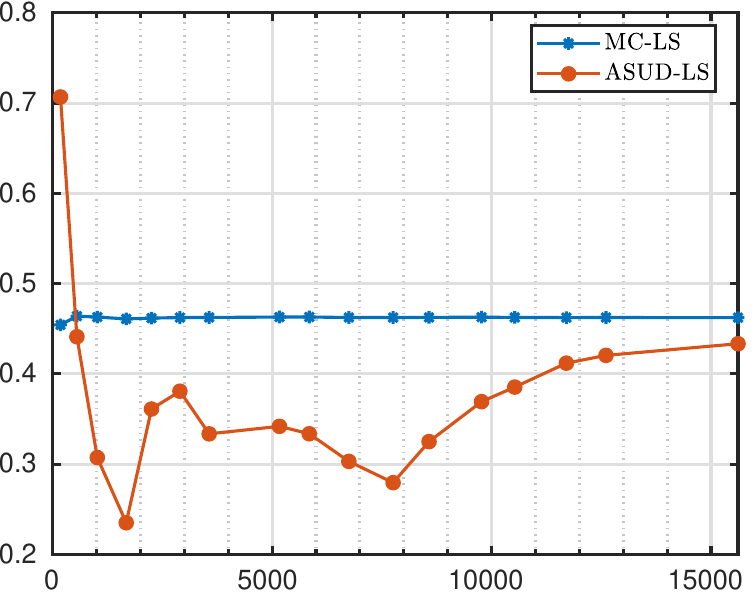} \\ 
\includegraphics[scale=0.35]{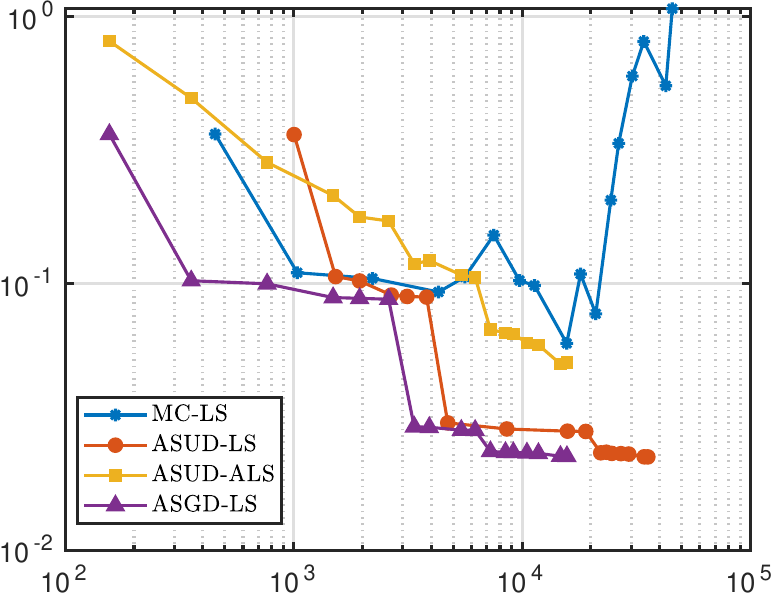} &  
\includegraphics[scale=0.35]{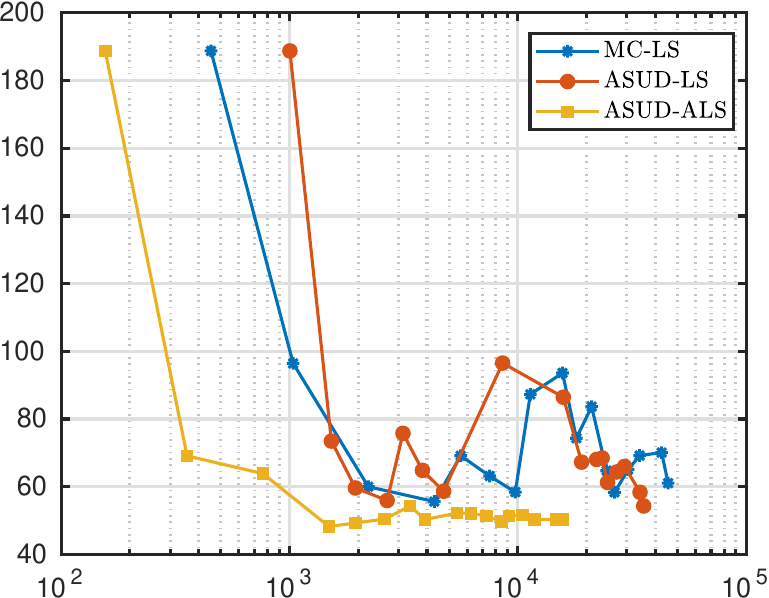} &  
\includegraphics[scale=0.35]{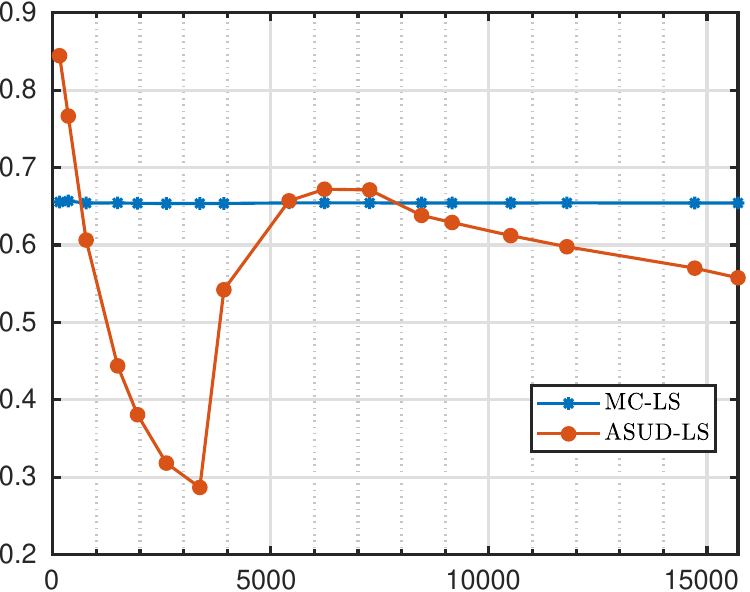} \\ 
\includegraphics[scale=0.35]{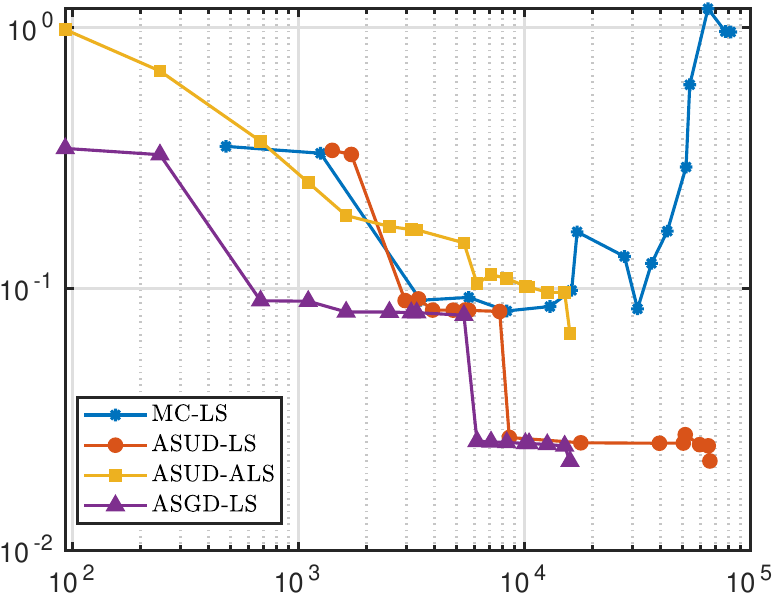} &  
\includegraphics[scale=0.35]{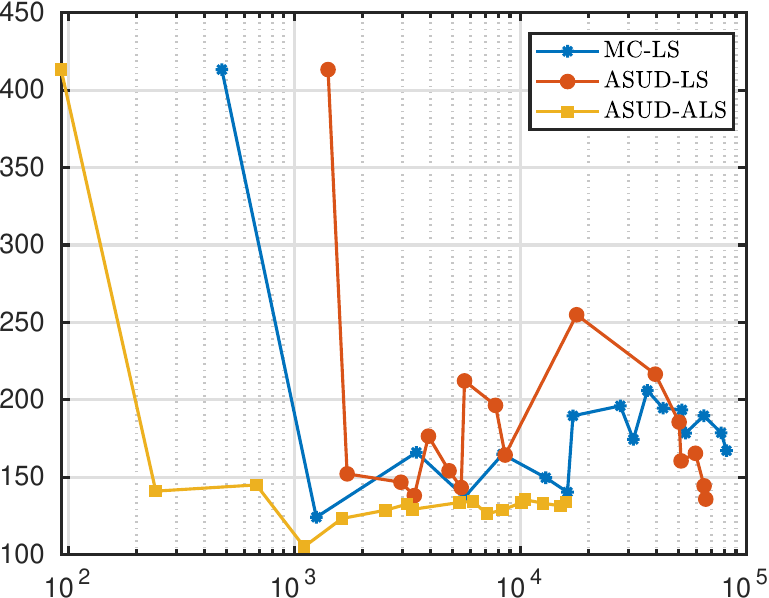} &  
\includegraphics[scale=0.35]{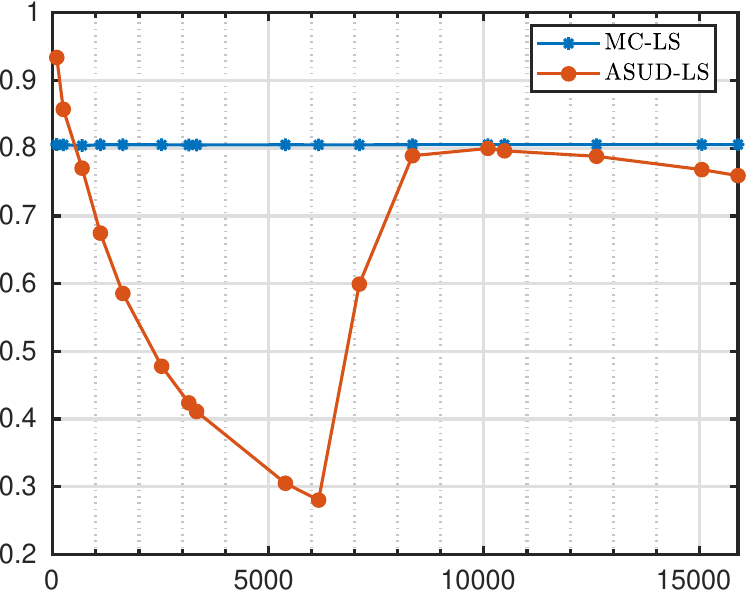}  
\end{tabular}
\end{center}
\vspace*{-1mm} 
\caption{Approximation of the function $f = f_2$ and domain $\Omega = \Omega_2$ in $d = 2,3,4,5,$ dimensions (top to bottom). Left: the relative error $E_l(f)$ versus the number of function evaluations $F_l$. Middle: the mismatch volume $V_l(f)$ versus $F_l$. Right: the rejection rate $R_l(f)$ versus $M_l$.
} 
\label{fig:F1_highdim}
\end{figure}

In Fig.\ \ref{fig:F1_highdim} we show experiments for $f = f_2$. Both ASUD-LS and ASUD-ALS can generally be seen to have decreasing error as we increase the number of functions evaluations used in computing their respective approximations, while this is not the case for MC-LS. This latter effect is indicative of a general property of Monte Carlo sampling with least-squares approximation: namely, it can be unstable and nonconvergent when the sampling rate is log-linear in the dimension of the polynomial subspace (note that in this case, $M_l \approx N_l \log(N_l)$, where $N_l = \dim(P_l)$). See \cite[Sec.\ 5.4.4]{adcock2021sparse} and \cite{MiglioratiThesis} for further discussion. Markedly, the errors of ASGD-LS and ASUD-LS follow a similar trend, with the latter requiring slightly more function evaluations due to lack of domain knowledge of $\Omega_2$ and the nonzero rejection rate (as shown in the right column of Fig.\ \ref{fig:F1_highdim}). Overall, ASUD-LS performs better than MC-LS when approximating the function, while both methods perform similarly when learning the domain. The mismatch volume is approximately $10\%$ for these methods in $d=2$ and $d=3$ dimensions. Notice also that for ASUD-ALS, the mismatch volume decreases faster than ASUD-LS and MC-LS in all dimensions studied. It is important to mention that this domain has decreasing volume relation to $D$ as the dimension increases. Consequently, the rejection rate for MC-LS increases with dimension as can be seen in the third column of Fig.\ \ref{fig:F1_highdim}. On the other hand, the rejection rate for ASUD-LS is generally lowe than Monte Carlo. This is clear in the lower-dimensional case $d=2$, but less noticeable in dimensions $d=3,4,5$.  

\begin{figure}[!h]
\begin{center}
\begin{tabular}{ccc} 
\includegraphics[scale=0.35]{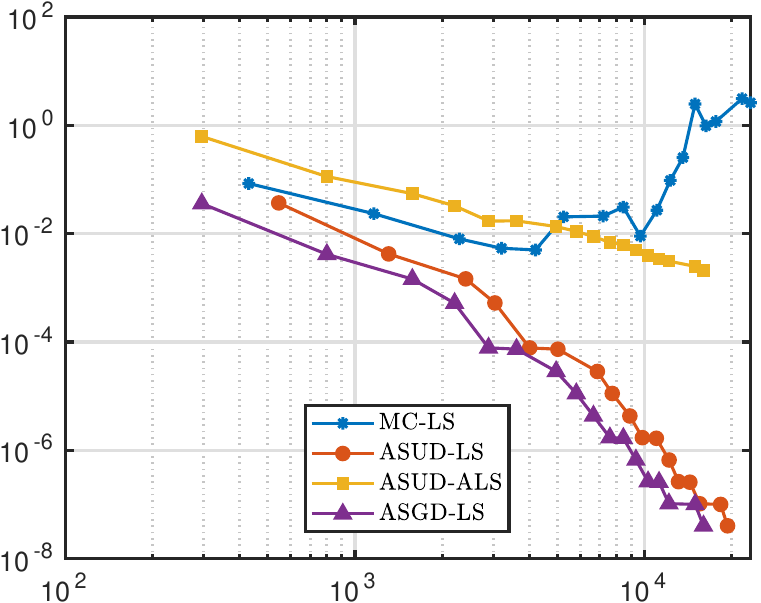} &  
\includegraphics[scale=0.35]{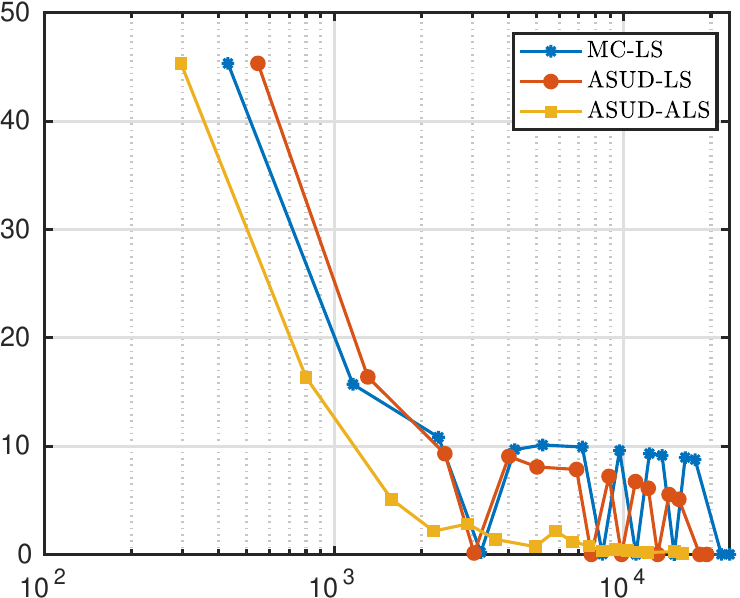} &  
\includegraphics[scale=0.35]{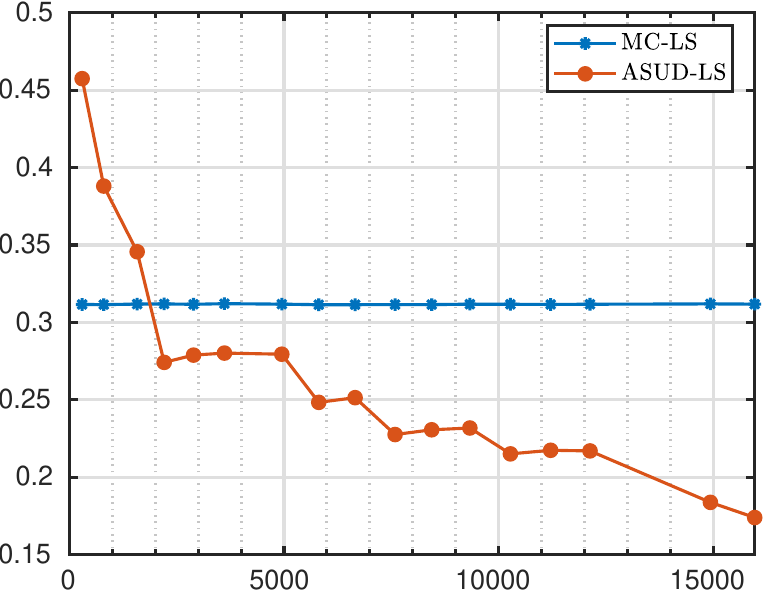} \\ 
\includegraphics[scale=0.35]{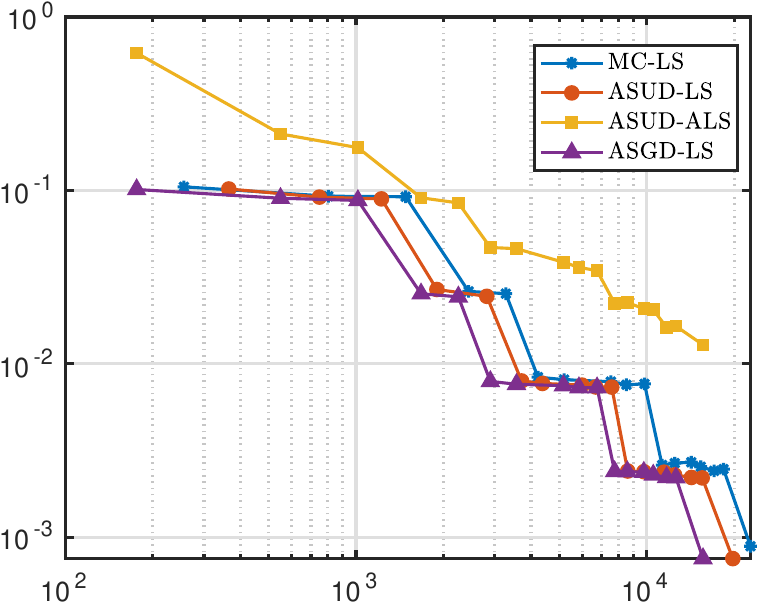} &  
\includegraphics[scale=0.35]{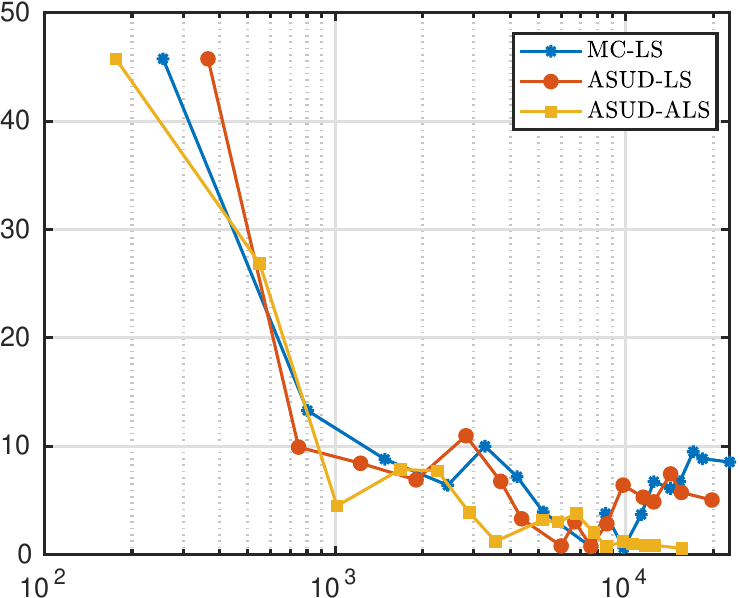} &  
\includegraphics[scale=0.35]{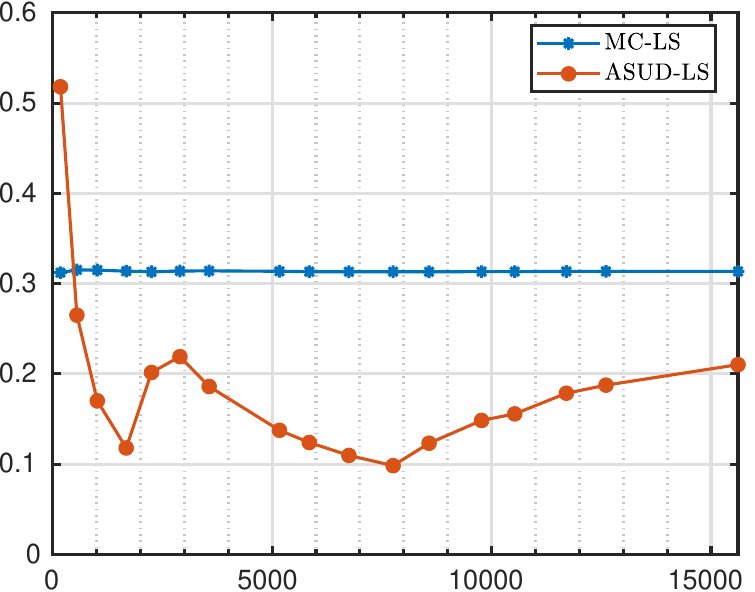} \\ 
\includegraphics[scale=0.35]{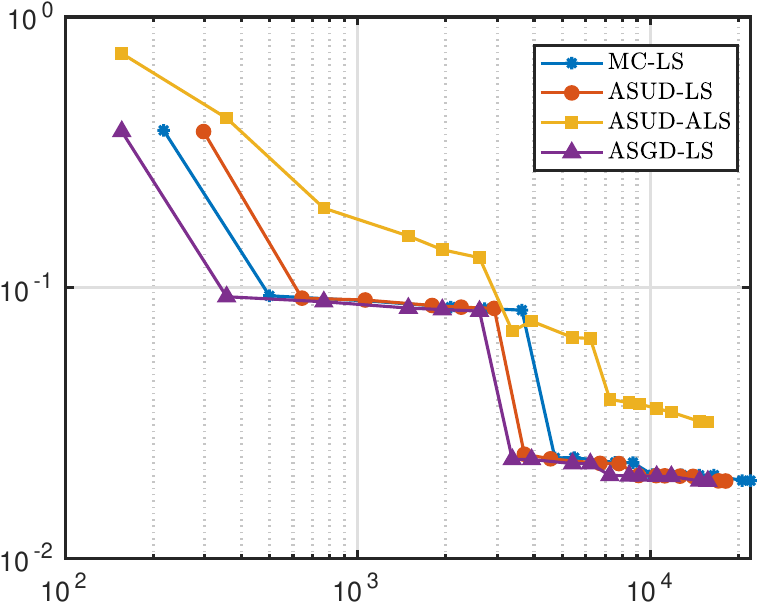} & 
\includegraphics[scale=0.35]{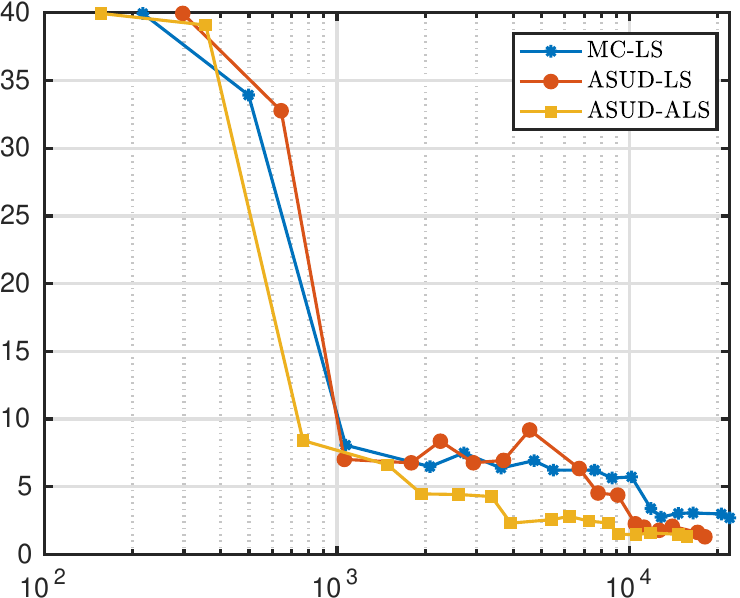} &  
\includegraphics[scale=0.35]{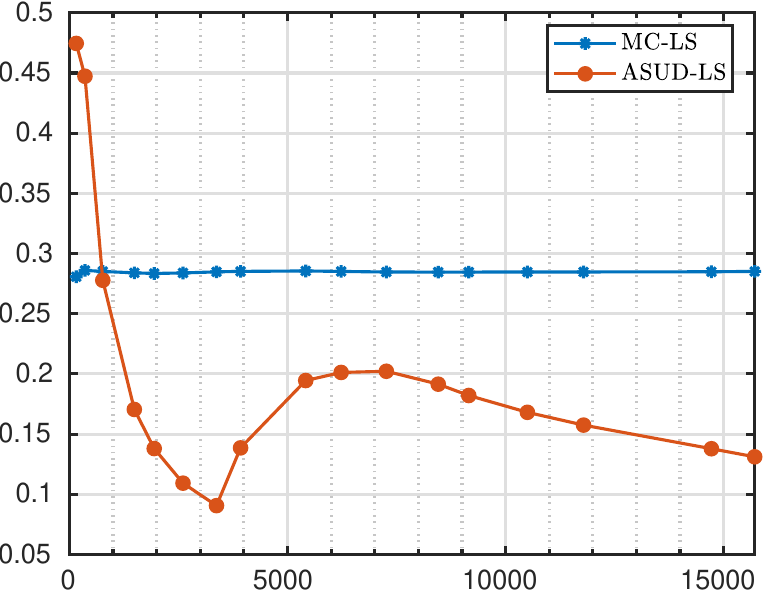} \\ 
\includegraphics[scale=0.35]{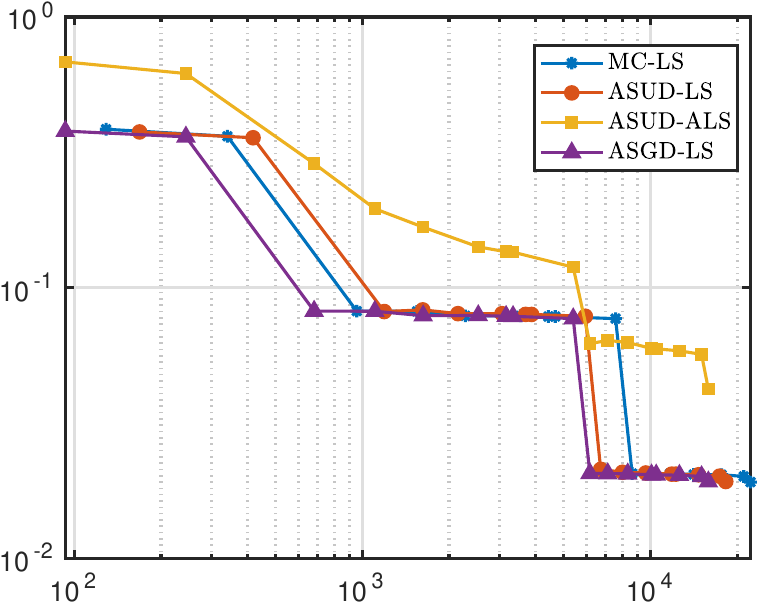} &  
\includegraphics[scale=0.35]{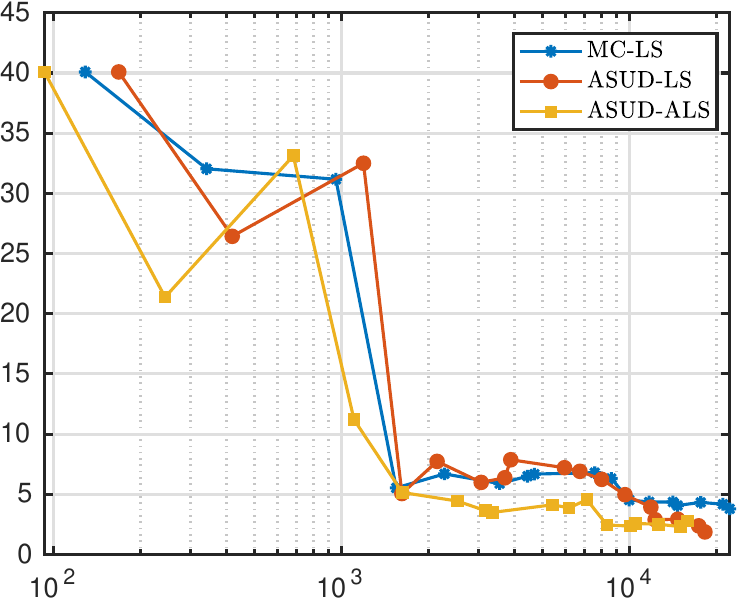} &  
\includegraphics[scale=0.35]{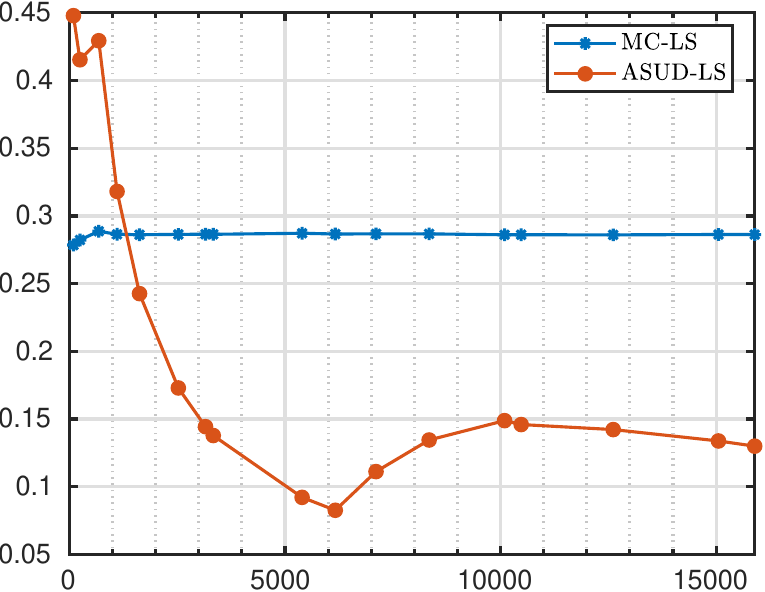} 
\end{tabular}
\end{center}
\vspace*{-1mm} 
\caption{Approximation of the function $f = f_3$ and domain $\Omega = \Omega_3$ in $d = 2,3,4,5,$ dimensions (top to bottom). Left: the relative error $E_l(f)$ versus the number of function evaluations $F_l$. Middle: the mismatch volume $V_l(f)$ versus $F_l$. Right: the rejection rate $R_l(f)$ versus $M_l$.
} 
\label{fig:F4_highdim}
\end{figure}

In Fig. \ref{fig:F4_highdim}, we show experiments for $f=f_3$. In two dimensions this function coincides with the function $f_2$. In higher dimensions, however, the volume of $\Omega_3$ relative to $D$ remains roughly constant, unlike in the case of $f_2$. Correspondingly, the rejection rate for MC-LS is roughly constant, and approximately 30\%. In each of the dimensions, ASUD-LS outperforms MC-LS for approximating the target function. As in the experiments in Fig.\ \ref{fig:F1_highdim} with $f=f_2$, the domain-learning performance is similar for the three methods, i.e.\ generally decreasing with increasing number of evaluations. However, ASUD-ALS decreases faster than both ASUD-LS and MC-LS.

\begin{figure}[!h]
\begin{center}
\begin{tabular}{ccc}  
\includegraphics[scale=0.35]{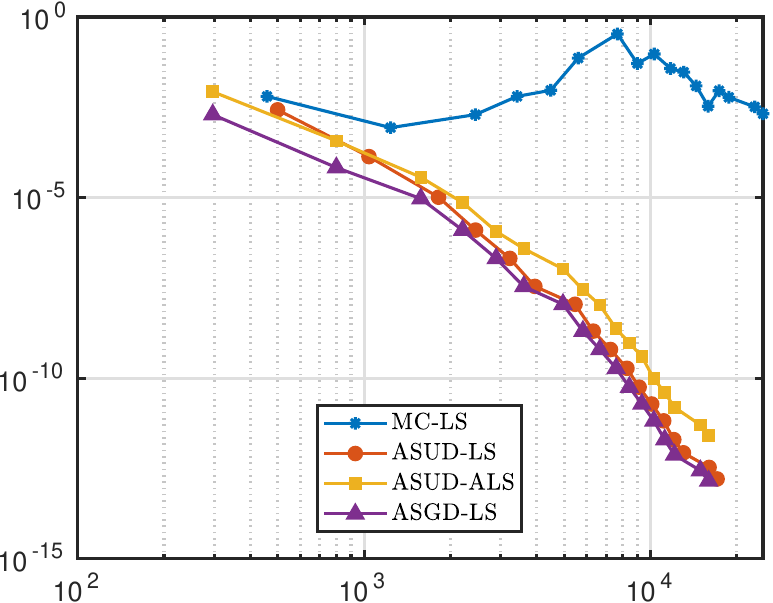} &  
\includegraphics[scale=0.35]{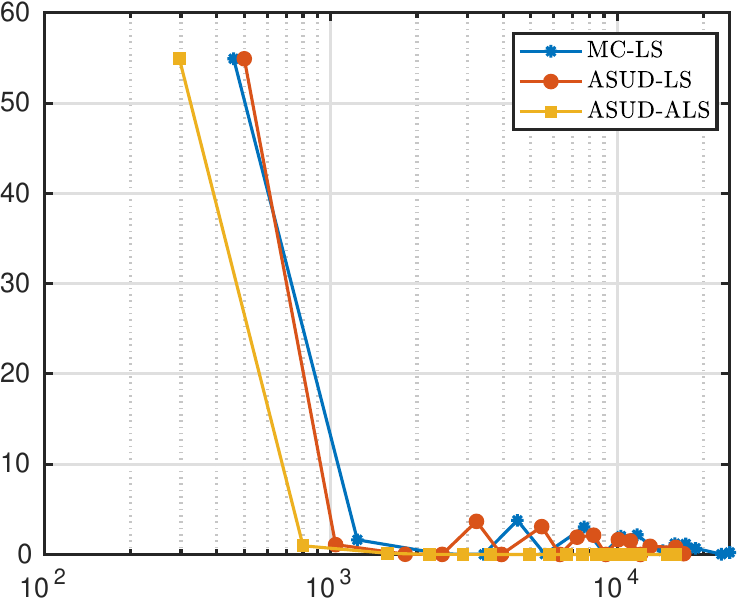} &  
\includegraphics[scale=0.35]{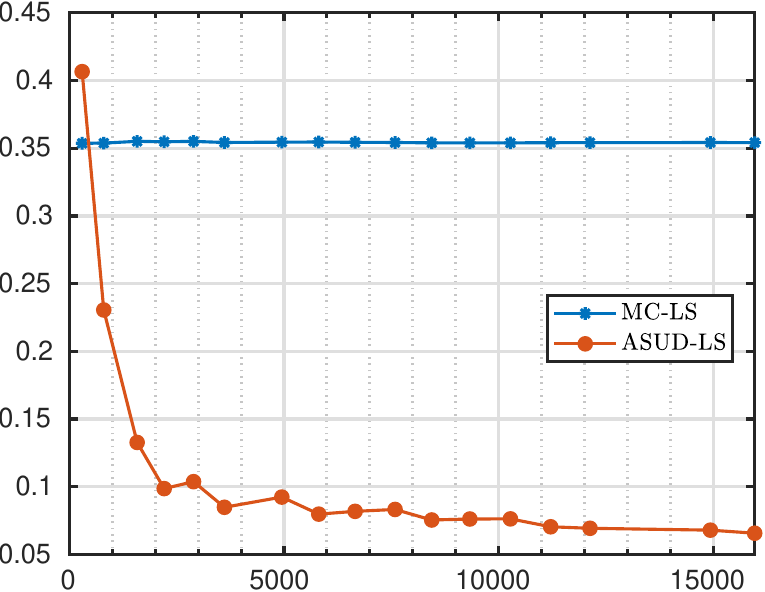} \\   
\includegraphics[scale=0.35]{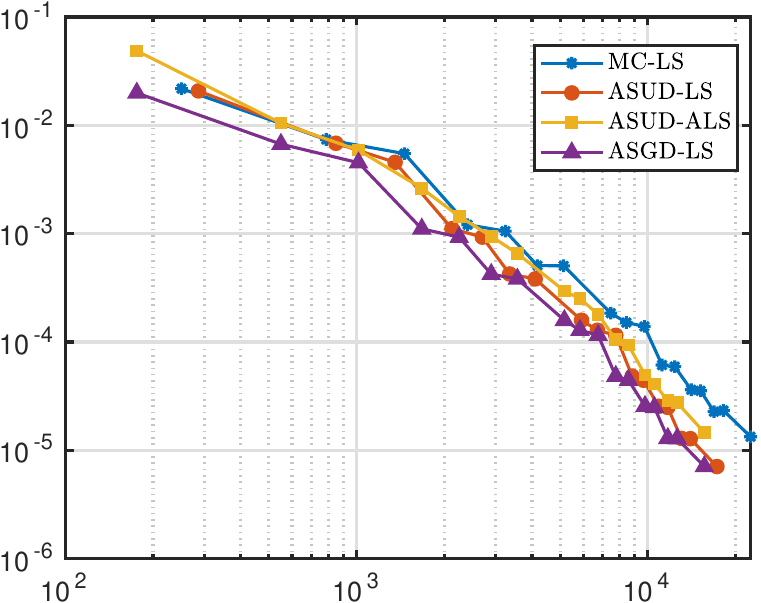} &  
\includegraphics[scale=0.35]{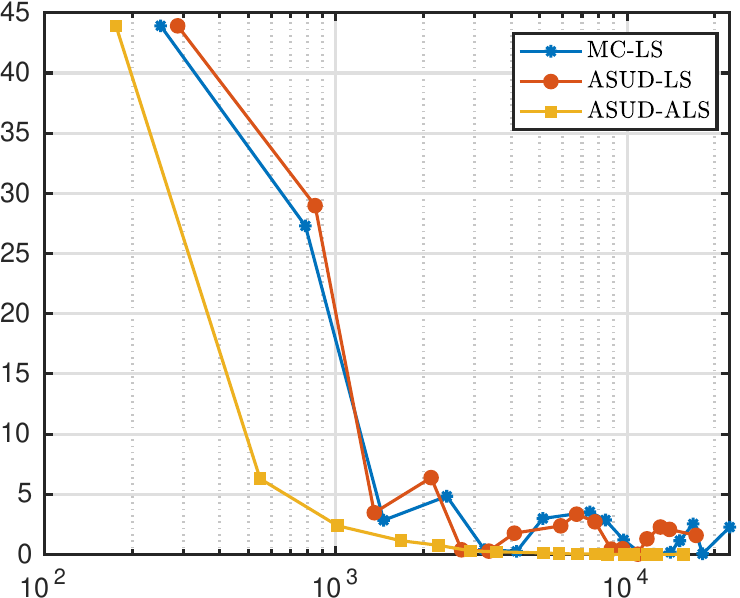} & 
\includegraphics[scale=0.35]{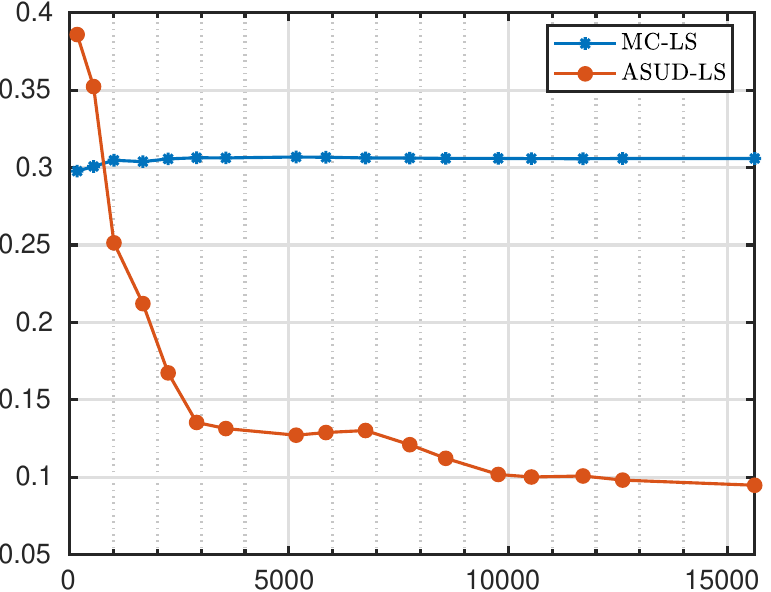} \\     
\includegraphics[scale=0.35]{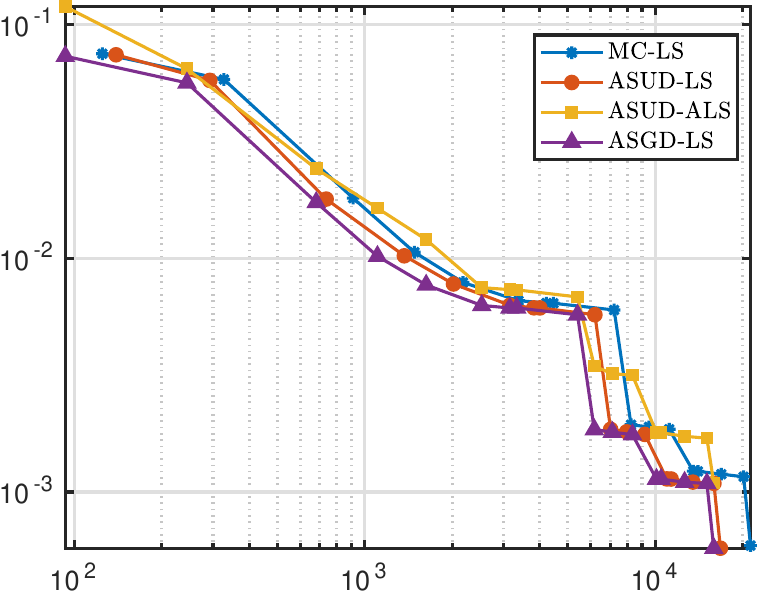} &  
\includegraphics[scale=0.35]{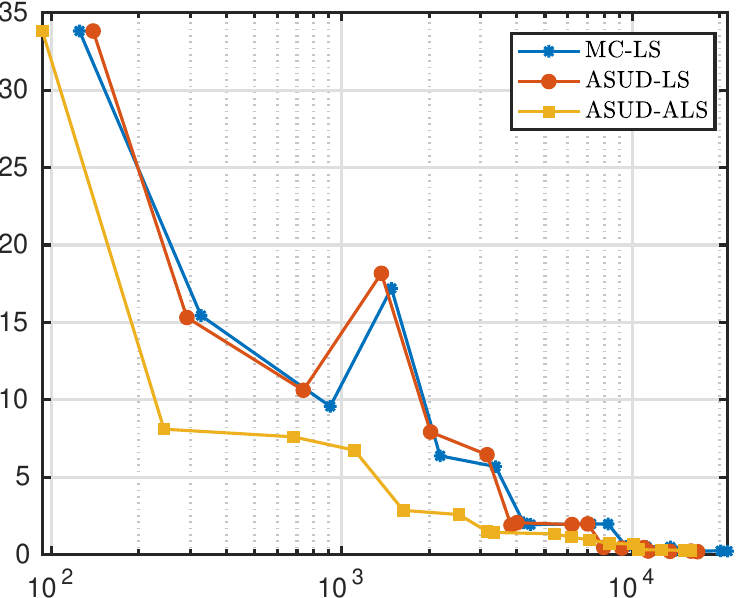} &  
\includegraphics[scale=0.35]{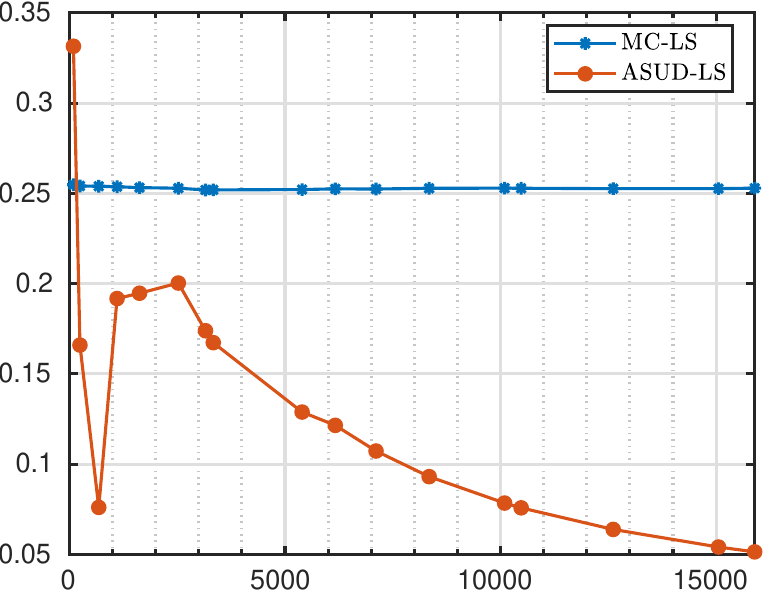} \\   
\includegraphics[scale=0.35]{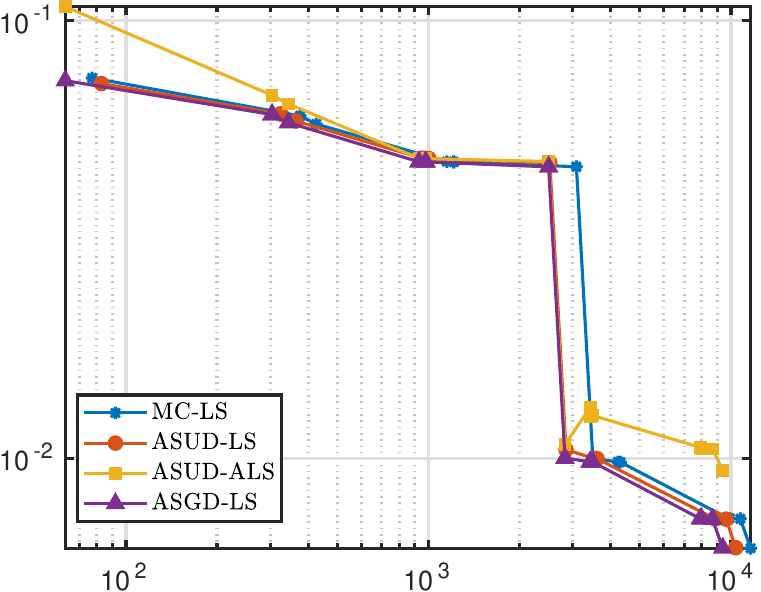} &  
\includegraphics[scale=0.35]{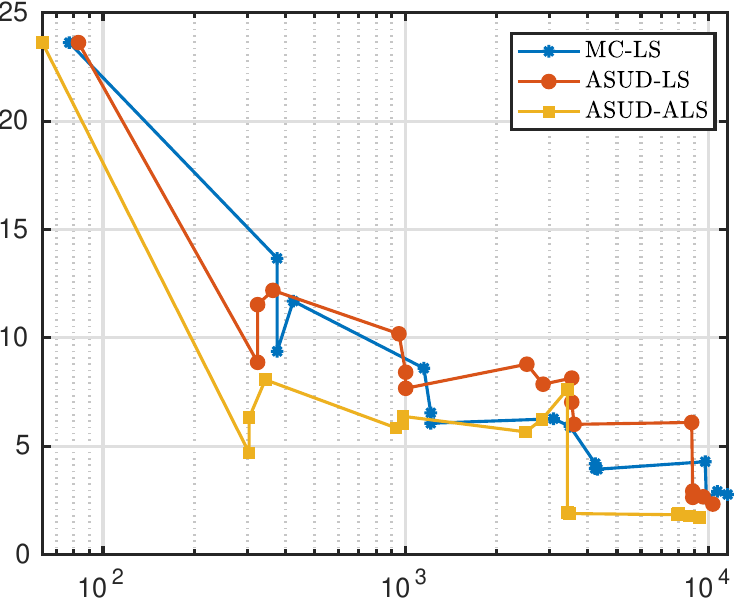} &  
\includegraphics[scale=0.35]{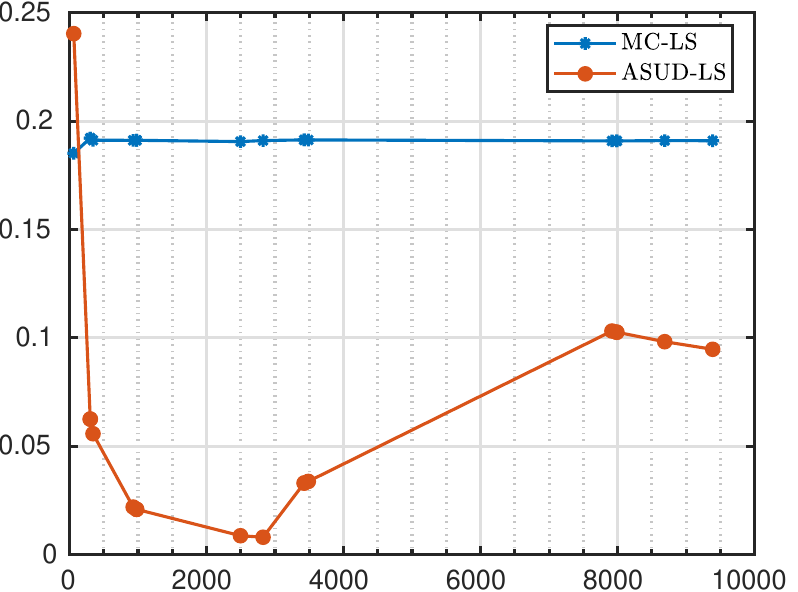} \\    
\includegraphics[scale=0.35]{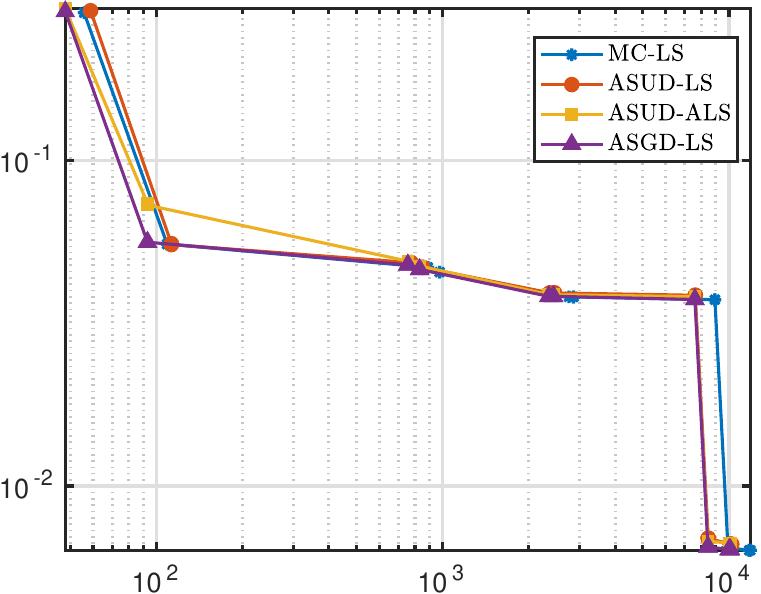} &  
\includegraphics[scale=0.35]{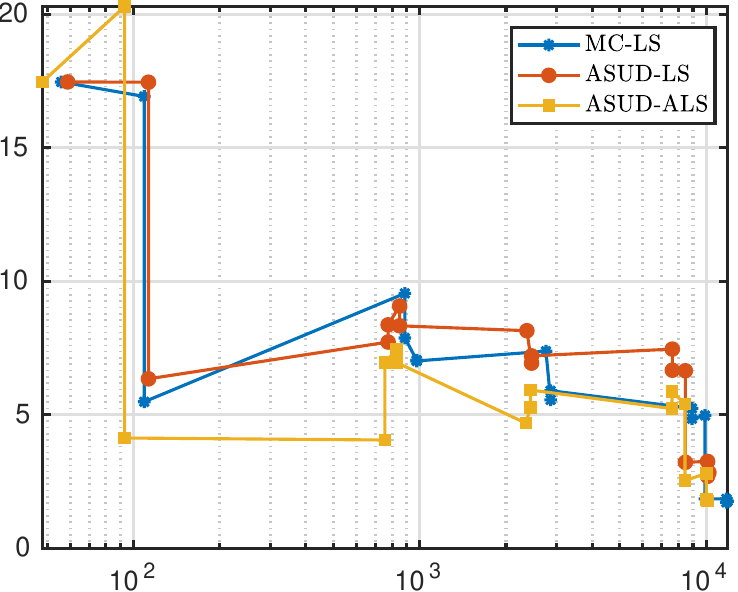} &  
\includegraphics[scale=0.35]{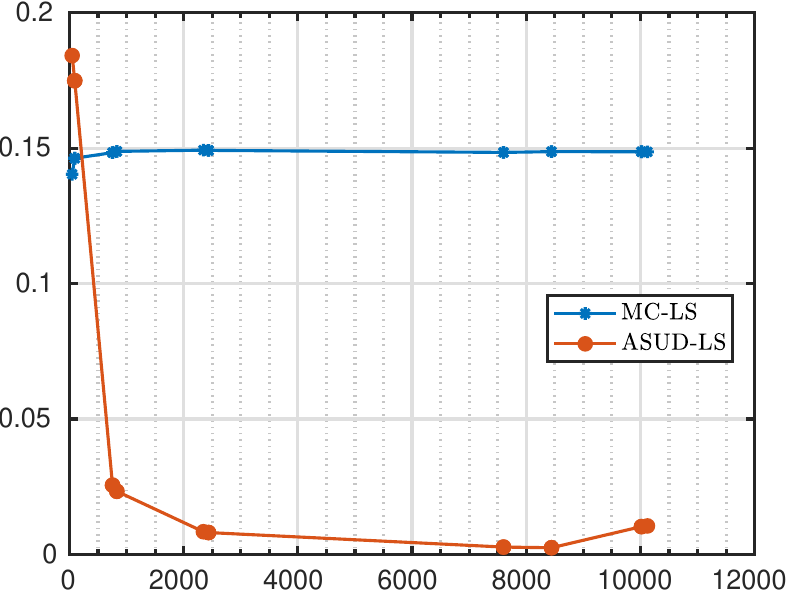} 
\end{tabular}
\end{center}
\vspace*{-1mm}
\caption{Approximation of the function $f = f_4$ and domain $\Omega = \Omega_4$ in $d = 2,3,5,10,15,$ dimensions (top to bottom). Left: the relative error $E_l(f)$ versus the number of function evaluations $F_l$. Middle: the mismatch volume $V_l(f)$ versus $F_l$. Right: the rejection rate $R_l(f)$ versus $M_l$.
} 
\label{fig:F6_highdim}
\end{figure}

Finally, in Fig.\ \ref{fig:F6_highdim} we show experiments for $f = f_4$. It is notable that ASUD-ALS performs similarly to ASUD-LS and ASGD-LS in terms of approximating the function in this case, and achieves better domain learning. We discuss the reasons behind this further in the next section.

\section{Theoretical discussion}\label{sec:theoretical}

\subsection{Accuracy and stability of weighted least squares}
 
Consider a domain $\Omega$ with a probability measure $\tau_{\Omega}$ and a weight function $w$ that is positive and finite almost everywhere on $\supp(\tau_{\Omega})$. Let $\bm{y}_1,\ldots,\bm{y}_M$ be $M$ sample points and consider the weighted least-squares approximation
\be{
\tilde{f} \in \argmin{p\in P}\left\lbrace \frac{1}{M}\sum_{i=1}^M w(\bm{y}_i) |f(\bm{y}_i) + e_i -p(\bm{y}_i))|^2\right\rbrace,
}
of a function $f\in L^2(\Omega,\tau_{\Omega})$ in a subspace $P \subseteq L^2(\Omega,\tau_{\Omega})$ with $\dim(P) = N \leq M$. Here, the values $e_i$ represent noise. Then the accuracy and stability (to noise) of this approximation are both controlled by the following \textit{discrete stability constant}
\bes{
\alpha = \inf \left \{ \nm{p}_{\mathrm{disc}} : p \in P,\ \nm{p}_{L^2(\Omega,\tau_{\Omega})} = 1 \right \},
}
where $\nm{g}_{\mathrm{disc}} = \sqrt{M^{-1} \sum^{M}_{i=1} w(\bm{y}_i) | g(\bm{y}_i) |^2 }$.
Specifically, if $\alpha > 0$ then the approximation is unique and it satisfies
\be{
\label{err-bd-LS}
\nmu{f - \tilde{f} }_{L^2(\Omega,\tau_{\Omega})} \leq \nmu{f - p}_{L^2(\Omega,\tau_{\Omega})} + \alpha^{-1} \nm{f - p}_{\mathrm{disc}} + \alpha^{-1} \tnm{\bm{e}},\quad \forall p \in P,
}
where $\tnm{\bm{e}} = \sqrt{\frac{1}{M} \sum^{M}_{i=1} w(\bm{y}_i) | e_i|^2}$ (see, e.g., \cite[Chpt.\ 5]{adcock2021sparse}). Thus, accuracy and stability of the approximation follow whenever $1/\alpha$ is not too large. Note that
\bes{
\nm{p}_{\mathrm{disc}} \geq \alpha \nm{p}_{L^2(\Omega,\tau_{\Omega})}.
}
Thus, when $1/\alpha$ is not too large, this states that the discrete norm over the sample points can be estimated from below by the underlying norm $\nm{\cdot}_{L^2(\Omega,\tau_{\Omega})}$ in which accuracy and stability of the approximation are estimated.

The discrete constant can be computed whenever the domain $\Omega$ is known. In that case, we compute an orthonormal basis $\{ \phi_1,\ldots,\phi_N \} \subseteq L^2(\Omega,\tau_{\Omega})$
of $P \subseteq L^2(\Omega,\tau_{\Omega})$ (if $\tau_{\Omega}$ is a finitely-supported measure, as it is in this paper, then we do this via QR decomposition, as described previously) and then use this to construct the matrix
\bes{
\bm{B} = \left ( \sqrt{w(\bm{y}_i)} \phi_j(\bm{y}_i) \right )^{M,N}_{i,j=1}.
}
The constant $\alpha$ is then precisely $\sigma_{\min}(\bm{B})$.

\begin{figure}[h!]
\begin{center}
{\small
\begin{tabular}{ccc} 
\includegraphics[scale=0.35]{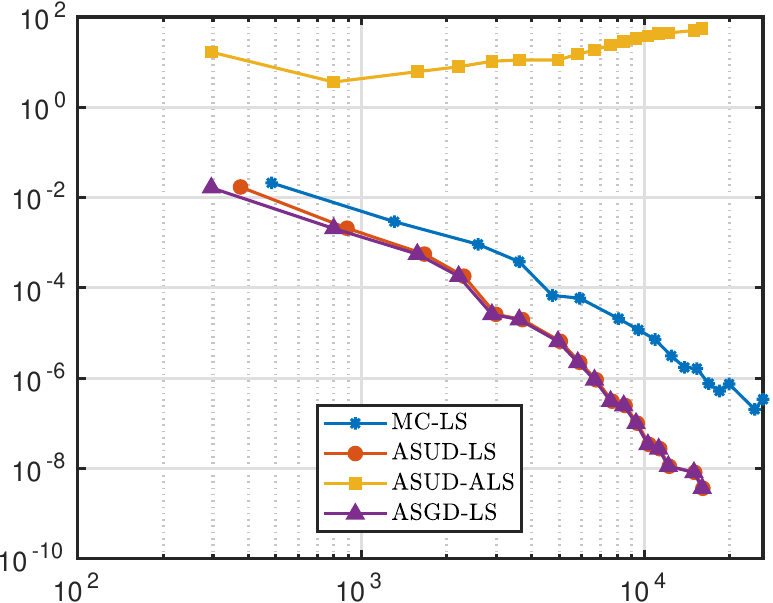} &  
\includegraphics[scale=0.35]{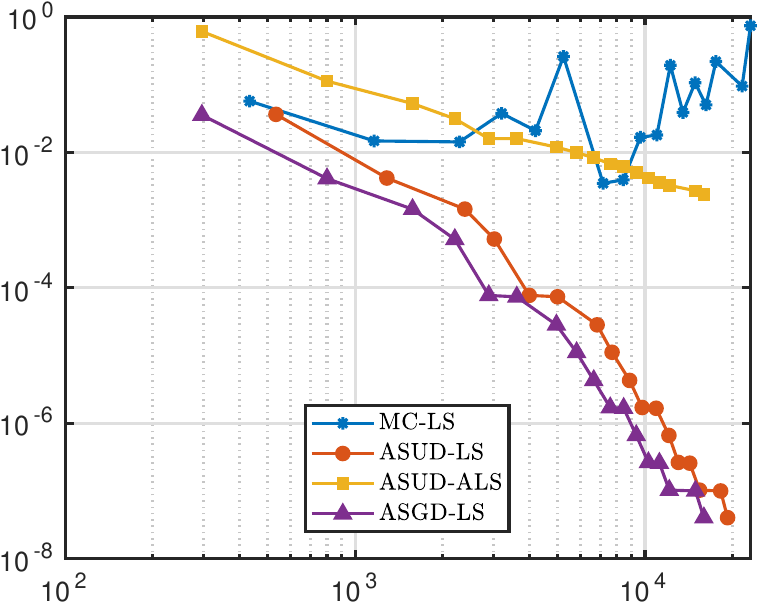} &  
\includegraphics[scale=0.35]{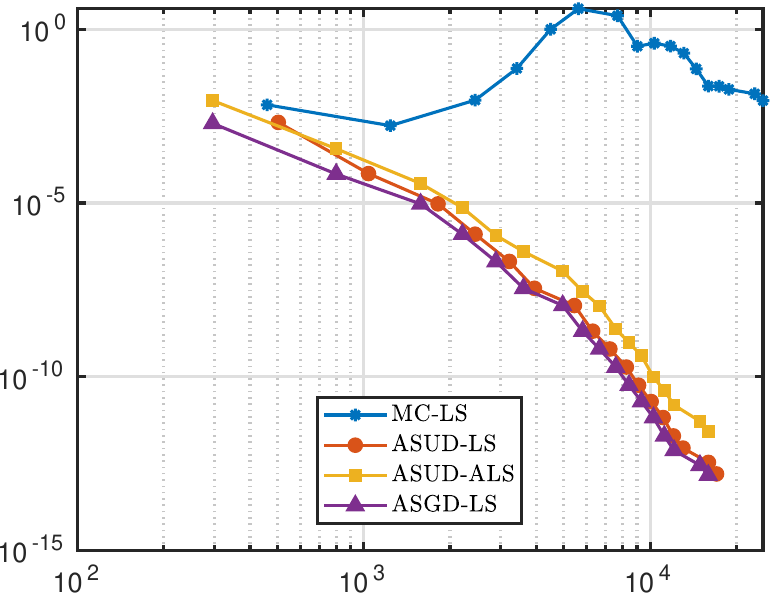} \\   
\includegraphics[scale=0.35]{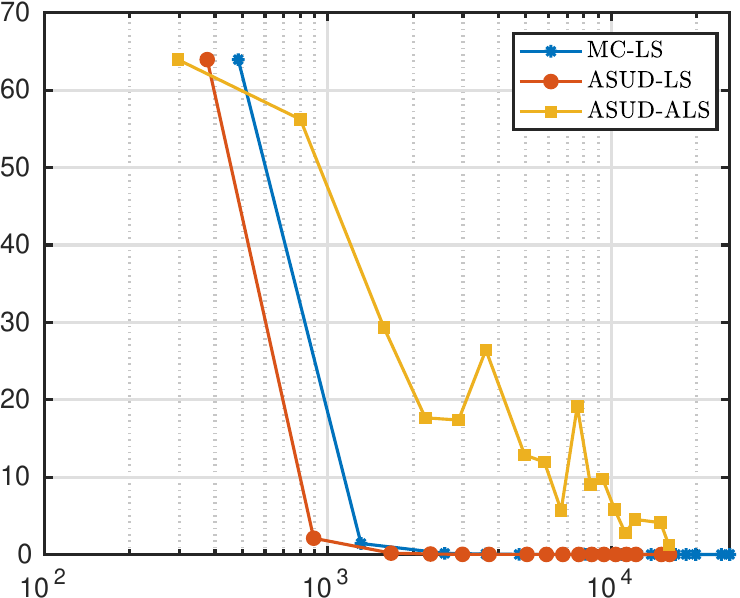} &  
\includegraphics[scale=0.35]{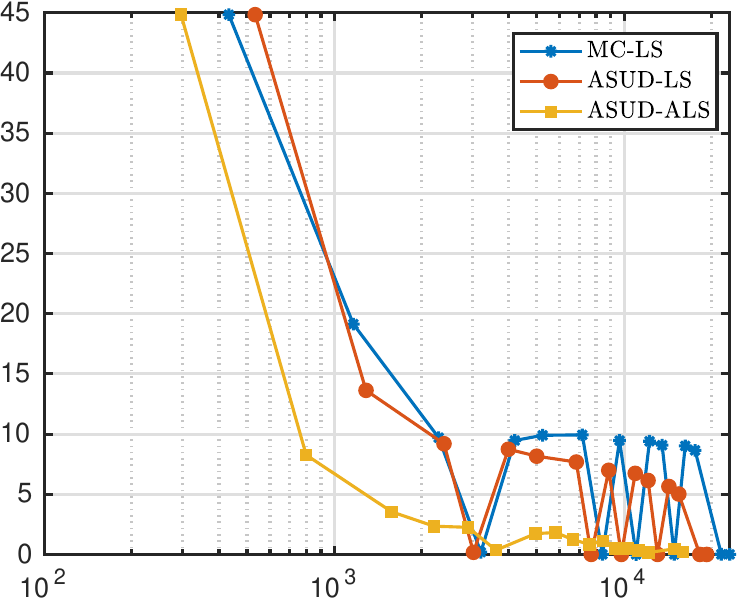} &  
\includegraphics[scale=0.35]{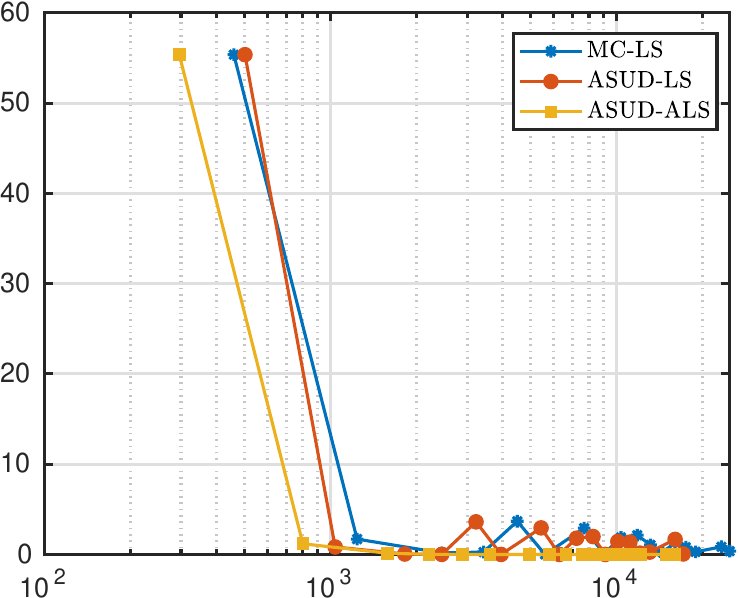} \\  
\includegraphics[scale=0.35]{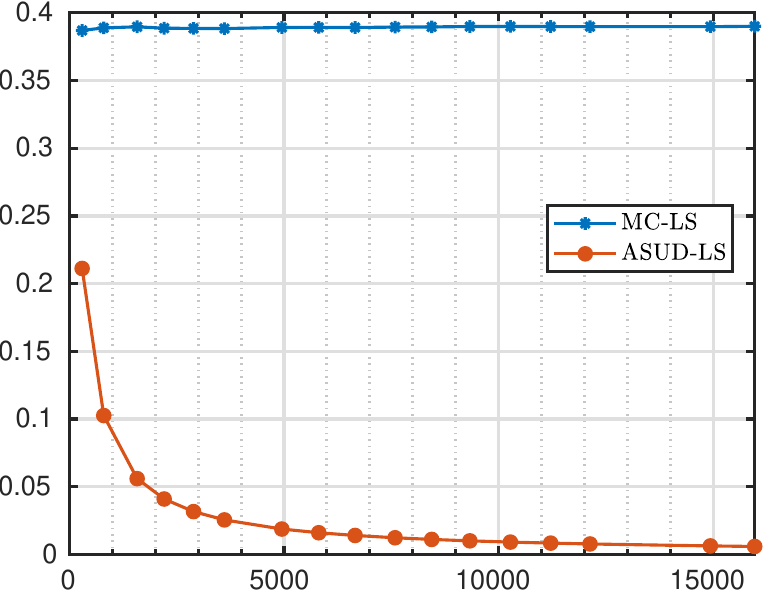} & 
\includegraphics[scale=0.35]{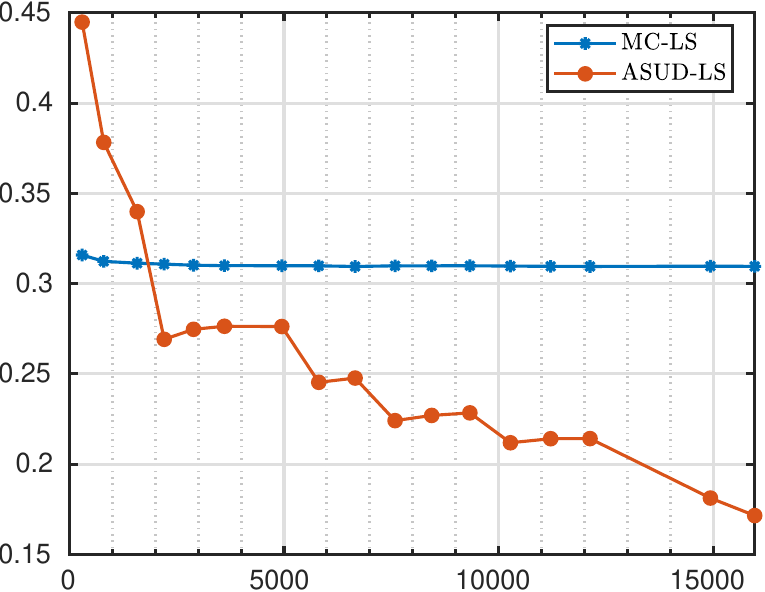} &  
\includegraphics[scale=0.35]{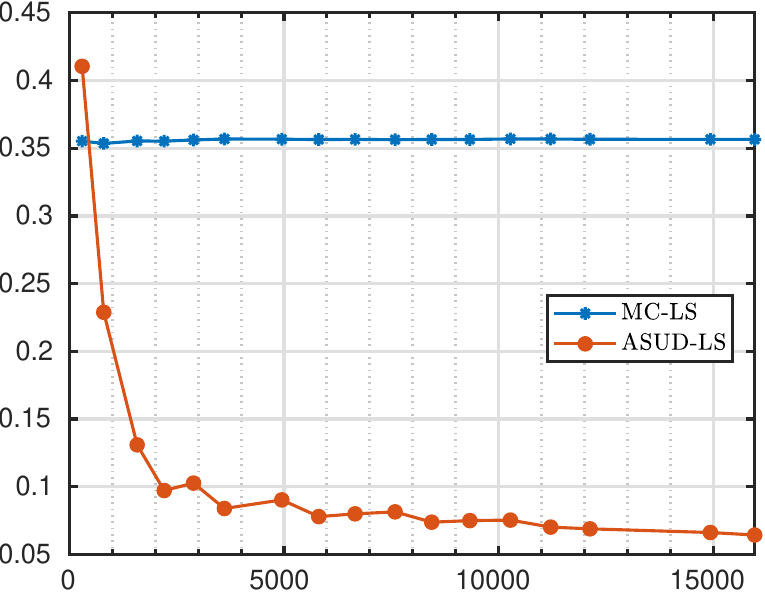} \\  
\includegraphics[scale=0.35]{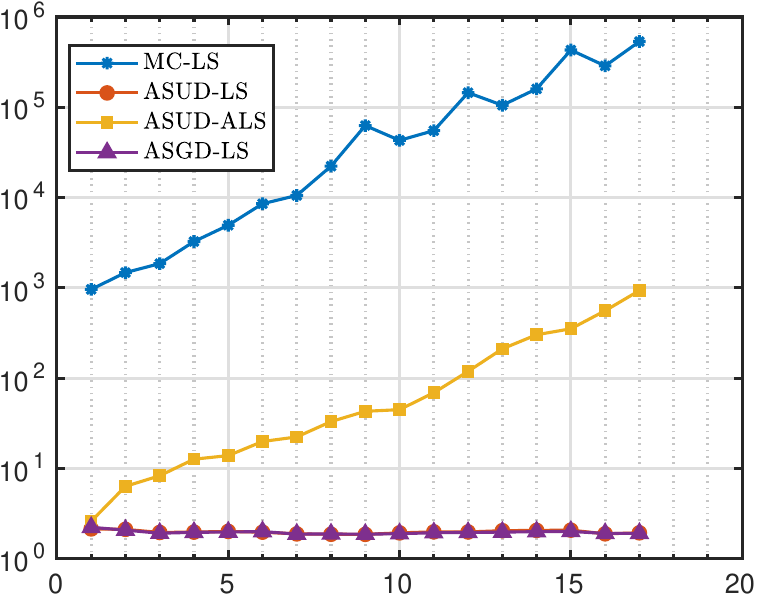} &  
\includegraphics[scale=0.35]{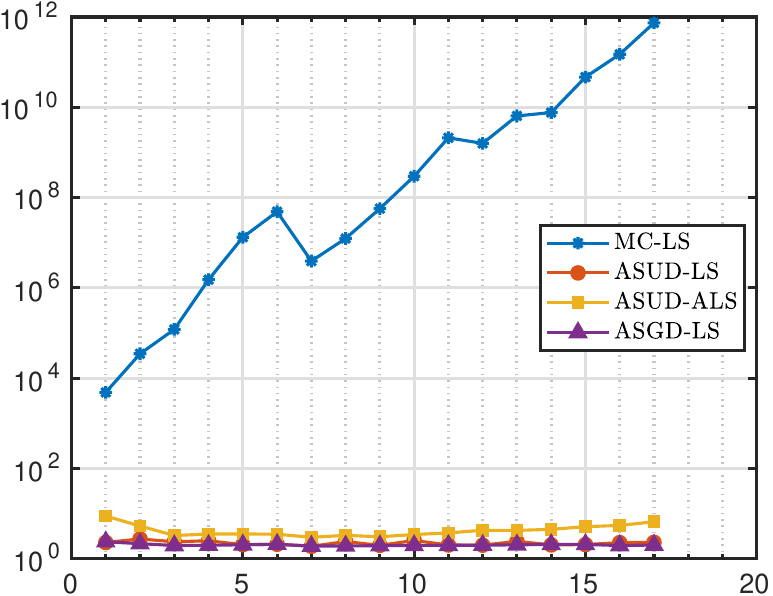} &  
\includegraphics[scale=0.35]{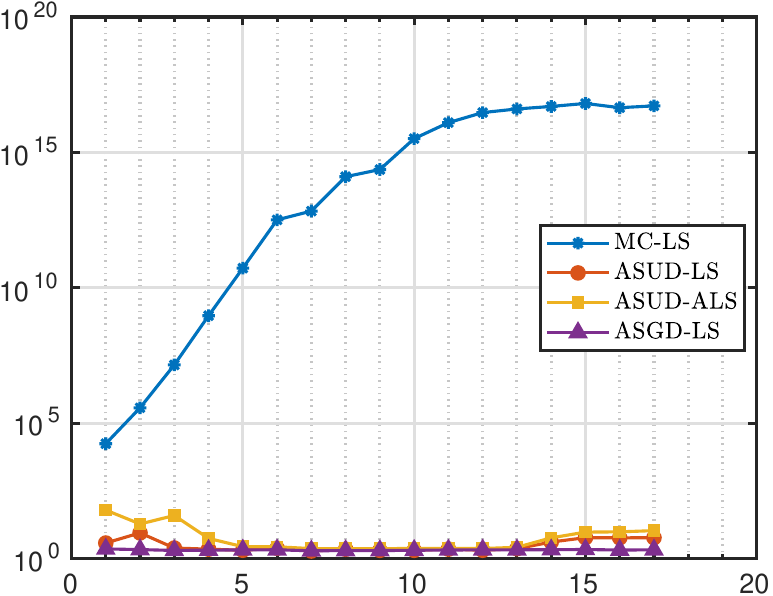} \\  
\includegraphics[scale=0.35]{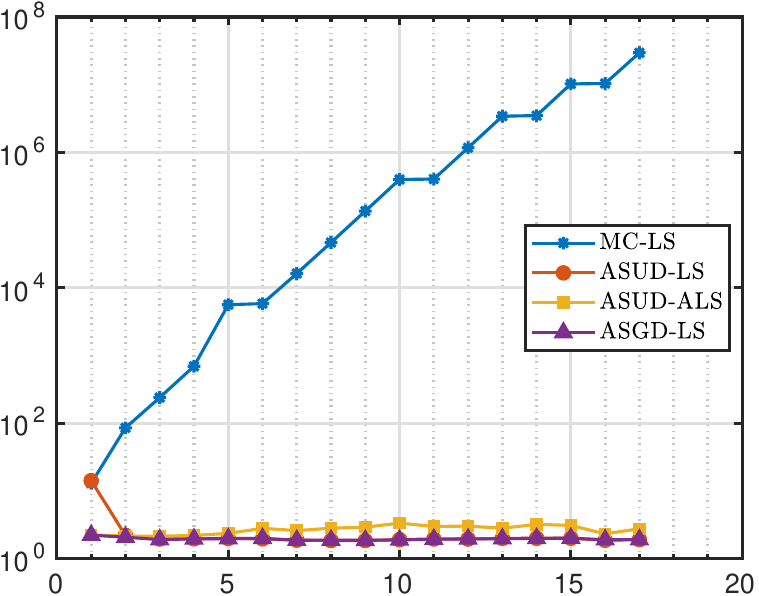} &  
\includegraphics[scale=0.35]{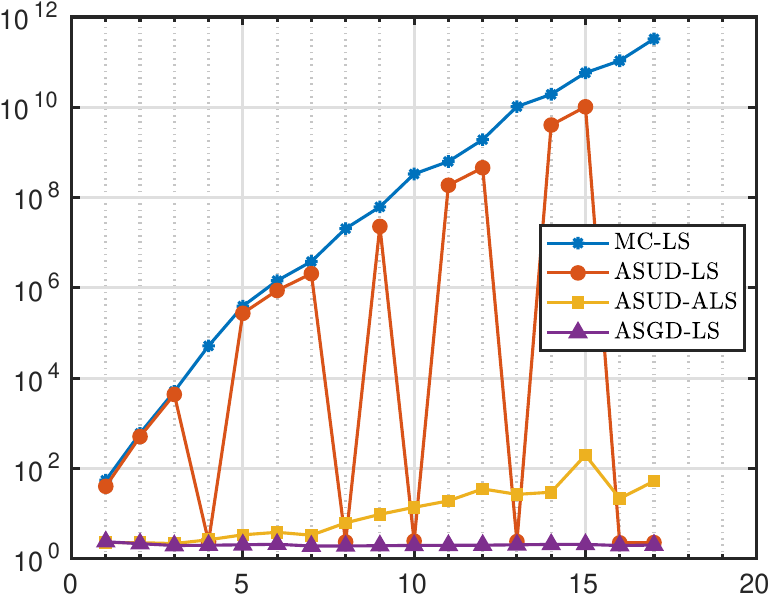} &  
\includegraphics[scale=0.35]{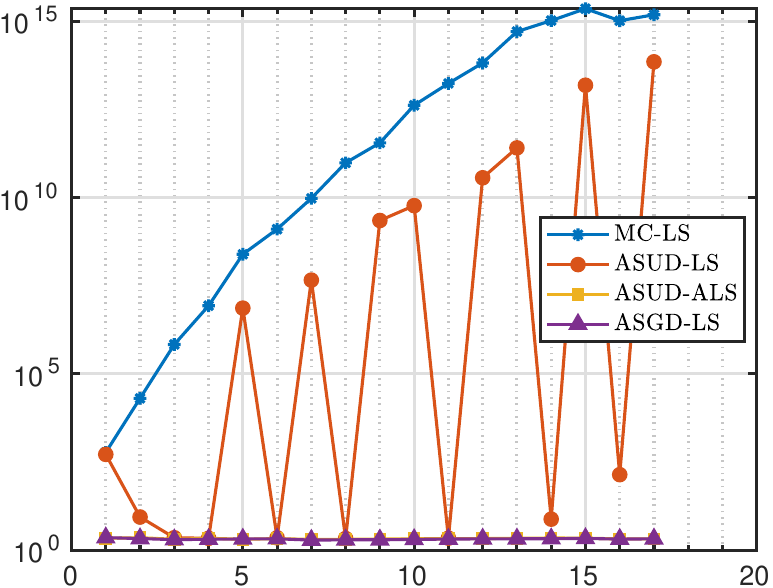}  
\end{tabular}
}
\end{center}
\vspace*{-1mm}
\caption{Approximation of the functions $f = f_1, f_2, f_4$ and domains $\Omega = \Omega_1, \Omega_2, \Omega_4$ (left to right) in $d = 2$ dimensions. First row: the relative error $E_l(f)$ versus the number of function evaluations $F_l$. Second row: the mismatch volume $V_l(f)$ versus $F_l$. Third row: the rejection rate $R_l(f)$ versus $M_l$. Fourth row: the constant $1/\alpha$ versus $l$. Fifth row: the constant $1/\beta$ versus $l$.
} 
\label{fig:WLS_AWLS_dim_2}
\end{figure}

In Fig.\ \ref{fig:WLS_AWLS_dim_2} we compute this constant for several different examples.  
As we see, in all cases MC-LS leads to a large constant that increases exponentially with the number of iterations (recall the earlier discussion). All three adaptive methods lead to much smaller constants, and therefore better stability. ASUD-LS  and ASGD-LS have similar constants, neither of which grow with the iteration number, and remain less than $10$ in magnitude. ASUD-ALS has a slightly large constant that can grow as large as roughly $10^2$ is size in these examples. 

This figure also shows the corresponding approximation errors. As we expect, the corresponding MC-LS approximation error is also larger than the ASUD-LS error, since the former has a much larger constant. It is notable that the ASUD-ALS error can be large, even though the constant $1/\alpha$ is small. The reason for this can be traced to the error bound \R{err-bd-LS}. For the MC-LS and ASUD-LS approximations, the term $\nm{f - p}_{\mathrm{disc}}$ is a discrete error over the sample points $\bm{y}_i$. Since in these methods the sample points are drawn from $\Omega$, this error is expected to be small whenever $f$ is well approximated over $\Omega$ by a function from $P$. However, in ASUD-ALS, the sample points are not restricted to belong to $\Omega$. They can, in theory, come from anywhere in $D$. Hence, the term $\nm{f - p}_{\mathrm{disc}}$ may be much larger in this case, since the functions considered are less smooth (in fact, in some cases, singular) over $D \backslash \Omega$. It is notable from Figs.\ \ref{fig:F1_highdim}, \ref{fig:F6_highdim} \& \ref{fig:F4_highdim} that ASUD-ALS gives the best relative performance in the task of domain learning for the functions $f_1$, $f_3$ and $f_4$, where $f_1$ and $f_4$ have logarithmic singularities at $\bm{y} = \bm{0}$ and $f_3$ is smooth on all of $D$. These results stand in contrast to those for function $f_2$, which has a quadratic singularity at $\bm{y}=\bm{0}$ and for which ASUD-ALS performs the worst, see Fig \ref{fig:F3_highdim}.

\subsection{Can accuracy and stability be numerically verified?}

As noted, the constant $\alpha$ can only be computed when the domain $\Omega$ is known. In this case, accuracy and stability of the approximation can numerically verified \textit{a priori}.  Unfortunately, this cannot be done when the domain is unknown, as it is in this paper. At a general step of ASUD-LS we compute a weighted least-squares approximation using an orthonormal basis $\{ \phi_1,\ldots,\phi_N \} \subseteq L^2(\Omega' , \tau')$, where $\Omega'$ an estimate for the true domain $\Omega$ and $\tau'$ is the normalized restriction of the measure $\tau$ on $D$ to $\Omega'$. If $\bm{y}_1,\ldots,\bm{y}_M$ are the sample points, then we form the least-squares matrix
\bes{
\bm{A} = \left ( \sqrt{w(\bm{y}_i)} \phi_j(\bm{y}_i) \right )^{M,N}_{i,j=1}.
}
It is clearly possible to compute the constant $\beta = \sigma_{\min}(\bm{A})$. Unfortunately, this may not describe the stability and accuracy of the approximation in the $L^2(\Omega,\tau_{\Omega})$-norm. Indeed, this constant is precisely
\bes{
\beta = \inf \left \{ \nm{p}_{\mathrm{disc}} : p \in P,\ \nm{p}_{L^2(\Omega',\tau')} = 1 \right \}.
}
Recall that the discrete norm $ \nm{p}_{\mathrm{disc}} $ is an approximation to the continuous norm $\nm{p}_{L^2(\Omega,\tau_{\Omega})}$ over $\Omega$. Thus, $\beta$ will only serve as a useful surrogate for the true stability constant $\alpha$ when $\Omega'$ is approximates the domain $\Omega$ sufficiently well. To see this, let $p \in P$. Then
\bes{
\nm{p}_{\mathrm{disc}} \geq \beta \nm{p}_{L^2(\Omega',\tau')}.
}
Suppose that $\Omega \subseteq \Omega'$, i.e.\ the true domain is a subset of the estimated domain. Then, if $c_{\Omega} = \int_{\Omega} \D \tau(\bm{y})$ and $c_{\Omega'} = \int_{\Omega'} \D \tau(\bm{y})$, we obtain the bound
\bes{
\alpha \geq \sqrt{c_{\Omega} / c_{\Omega'}} \beta \geq \sqrt{c_{\Omega}} \beta,
}
where in the second step we use that fact that $c_{\Omega'} \leq 1$, since $\Omega' \subseteq D$. Hence, in this case, we expect $\beta$ to provide a reasonable surrogate for $\alpha$. Unfortunately, if $\Omega \not \subseteq \Omega'$ then this is not the case. Indeed, in general we have
\bes{
\alpha \geq   \beta \gamma,\qquad \gamma : =  \inf \left \{ \nm{p}_{L^2(\Omega',\tau')} : p \in P, \ \nm{p}_{L^2(\Omega,\tau_{\Omega})} = 1 \right \}.
} 
The latter term can easily be large unless $\Omega'$ is a very good estimate of $\Omega$. For simplicity, consider the case $\Omega' \subseteq \Omega$. Then $1/\gamma$ determines how large an element $p \in P$ can grow on the larger domain $\Omega$ in relation to its size over $\Omega'$. In the case of polynomial subspaces, this is essentially a type of Remez inequality. Such inequalities are known to grow exponentially in the polynomial degree with a rate depending on the relative volumes of the two domains  (see, e.g., \cite{TemlyakovRemezHC,GanzburgRemez}).

In Fig.\ \ref{fig:WLS_AWLS_dim_2} we also show the constant $1/\beta$ for the various approximations. For MC-LS it grows large, much like the constant $1/\alpha$, since the approximation is unstable. For ASUD-LS the constant can also grow large, depending on the problem and the domain estimate. Indeed, comparing it with the domain mismatch volume, we see that the size of the behaviour of this constant closely tracks with when the domain estimate gets worse. On the other hand, $1/\beta$ remains small for the ASUD-ALS approximation. This is indicative of the fact that this scheme performs an approximation over the whole estimated domain. However, as we have seen, this scheme often leads to worse approximation errors for precisely this reason.

The main conclusion of this section is that computing the constant $\beta$ is of limited value for ASUD-LS, since $1/\beta$ being large does not necessarily imply a poor approximation (to either the function or the domain) and $1/\beta$ being small need not imply a good approximation (to either the function or the domain). 
 
\section{Conclusions}
\label{sec:conclusions}

In this paper, we introduced a new method, ASUD, for function approximation and domain learning over unknown domains of interest. This method combines previous work on weighted least-squares approximation on general domains (ASGD) with a domain estimation procedure. As shown in our numerical experiments, this procedure can lead to significant advantages over standard Monte Carlo sampling (MC-LS), even in higher dimensions.

Since it employs rejection sampling, a limitation of this approach is that it may suffer a high rejection rate in the first step if the domain $\Omega$ has very small measure in comparison to $D$. However, it is unclear whether or not this problem can be avoided at the level of generality considered in the paper. On the other hand, if an initial estimate of $\Omega$ is known, then this can be used to reduce the initial rejection rate by making a better informed choice of $D$.

There are several avenues for further investigations. A first one is to enhance the sampling efficiency further. In recent work \cite{haberstich2019boosted}, a boosting procedure for weighted least-squares approximations has been introduced. This procedure uses resampling and a greedy strategy to selectively remove sample points. Incorporating this into ASUD has the potential to further enhance its performance. Another topic for future research involves changing the approximation scheme. In this work we have used (weighted) least squares based on a fixed sequence of subspaces $P_1 \subseteq P_2 \subseteq \cdots$. An interesting extension, as mentioned briefly in \S \ref{sec:ASUD}, is the case where the subspaces are also generated adaptively, via, for instance, greedy adaptive methods \cite{cohen2018multivariate,migliorati2015adaptive,migliorati2019adaptive}. A related approach involves using compressed sensing via (weighted) $\ell^1$-minimization (see, e.g., \cite{adcock2018compressed,adcock2021sparse}). See \cite{adcock2022towards} for recent work on polynomial approximation via compressed sensing on irregular domains. Finally, there is also increasing interest in using deep learning for high-dimensional approximation tasks arising in UQ (see, e.g., \cite{adcock2021deep,adcock2021gap} and references therein). An interesting open problem is to combine the adaptive sampling procedure developed in this paper with a suitable deep neural network training strategy. See \cite{AdcockEtAl2022} for a recent adaptive sampling strategy based on Christoffel functions for deep learning on known domains.

\section*{Acknowledgements}

The authors acknowledges support from the Natural Sciences and Engineering Research Council of Canada (NSERC) through grant 611675 and the Pacific Institute for the Mathematical Sciences (PIMS).

%%%%%%%%%%%%%%%%%%%%%%%%%%%%%%%%%%%%%%%%%%%%%%%%%%%%%%%%%%%%%%%%%%%%%%%
%%%----------------        References     		    ----------------%%%
%%%%%%%%%%%%%%%%%%%%%%%%%%%%%%%%%%%%%%%%%%%%%%%%%%%%%%%%%%%%%%%%%%%%%%%

\bibliographystyle{abbrv}
\small
\bibliography{Adaptive_SamplingRefs}

\end{document}